\documentclass[11pt]{article}
\usepackage[T1]{fontenc}
\usepackage{lmodern}
\usepackage{amsmath,amssymb,amsthm}
\usepackage[margin=1in]{geometry}
\usepackage{hyperref}
\usepackage{float}
\newtheorem{theorem}{Theorem}[section]
\newtheorem{proposition}[theorem]{Proposition}

\newtheorem{lemma}[theorem]{Lemma}

\theoremstyle{remark}
\newtheorem{remark}[theorem]{Remark}
\theoremstyle{plain}
\newtheorem{corollary}[theorem]{Corollary}

\usepackage{booktabs,color}
\usepackage{mathrsfs}
\title{The Classification of Orthogonal Polynomials on the Unit Circle
Satisfying a Linear Second-Order Differential Equation with
Varying Polynomial Coefficients}

\author{Jorge A. Borrego--Morell\\
\small Universidade Federal do Rio de Janeiro\\
\small Campus Duque de Caxias, Brazil\\
\small ORCID: \href{https://orcid.org/0000-0001-8871-2739}
{0000-0001-8871-2739}}

\date{}

\begin{document}

\maketitle

\begin{abstract}
We establish a Bochner-type classification for positive
probability measures of infinite support on the unit circle
whose monic orthogonal polynomials satisfy
\[
    A(z)\Phi_n''(z)+B(z,n)\Phi_n'(z)+C_n\Phi_n(z)=0,
    \qquad n\geq1,
\]
where $A\not\equiv0$ is fixed, $\deg A\leq2$,
$\deg_z B\leq1$, and $C_n\in\mathbb C$.
The classification is carried out for the canonical
configurations of the leading coefficient
$\theta(z)\in\{1,z,z(1-z),z^2\}$ and, up to a common
invertible complex affine change of variable, for their
affine equivalents, with the transformation acting
simultaneously on the differential equation and the
orthogonality measure.
Our approach relies on coefficient comparison, the location
of the zeros of the orthogonal polynomials, and compatibility
with the Szeg\H{o} recurrence. 
\end{abstract}



\section{Introduction}
\label{sec:introduction}

The characterization of orthogonal polynomials by differential
equations is a central problem in the theory of special functions.
In his foundational work, Bochner~\cite{Bochner1929} classified
polynomial eigensystems of a common second-order differential
operator with polynomial coefficients. Up to complex affine
changes of the independent variable, his classification comprises
the Hermite, Laguerre, Jacobi, and generalized Bessel families,
together with the exceptional monomial family; see
also~\cite[pp.~974--975]{KwonLittlejohn1997}.

Under the additional requirement of orthogonality with respect to
a positive measure on the real line with infinite support, only
the Hermite, Laguerre, and Jacobi families remain, subject to the
usual parameter restrictions. When positivity is relaxed and
quasi-definite moment functionals are considered, the generalized
Bessel family also occurs, together with broader parameter ranges
of the Laguerre and Jacobi families~\cite{KwonLittlejohn1997}.
In particular, real sign-changing orthogonalizing weights for
generalized Bessel polynomials were constructed
in~\cite{EvansEtAl1993}.

In this paper we consider positive orthogonality measures whose monic
polynomials satisfy
\begin{equation}
\label{eq:intro-degree-dependent-equation}
    A(z)\Phi_n''(z)+B(z,n)\Phi_n'(z)+C_n\Phi_n(z)=0,
    \qquad n\geq1,
\end{equation}
where $A\not\equiv0$ is independent of $n$, $\deg A\leq2$,
$\deg_z B\leq1$, and $C_n\in\mathbb C$. If
$A(z)=a_2z^2+a_1z+a_0$ and $B(z,n)=\beta_{1,n}z+\beta_{0,n}$,
comparison of leading terms gives
\[
    C_n=-n\bigl((n-1)a_2+\beta_{1,n}\bigr).
\]
The problem is to determine which polynomial equations of this form
are compatible with a complete positive OPUC sequence and to
identify the corresponding measures. We refer to this formulation
as a \emph{Bochner-type problem on the unit circle}.

The degree dependence of differential coefficients has a classical
precedent in the work of Shohat~\cite{Shohat1939} on what are now
called semiclassical orthogonal polynomials: their second-order
differential equations may have degree-dependent leading
coefficients. In the spirit of Bochner's original
characterization, the independence of $A$ from $n$ identifies a
common principal differential structure for the entire polynomial
sequence, thereby specifying the notion of classicality adopted
in the present work. The degree dependence of $B_n$, on the other
hand, allows nontrivial positive families on the unit circle.

Several approaches to differential equations for OPUC have been
developed. Ismail and Witte~\cite{IsmailWitte2001} constructed
raising and lowering operators, derived second-order differential
equations, and obtained functional identities, discriminant
formulas, and an electrostatic interpretation of the zeros. These
constructions provide general tools for deriving differential
relations from a sufficiently regular weight. Branquinho and
Rebocho~\cite{BranquinhoRebocho2009} characterized Laguerre--Hahn
sequences on the unit circle through second-order matrix
differential equations associated with a Riccati equation for the
Carath\'eodory function. This relates differential properties of the
measure itself to differential systems for the polynomials and
their associated functions.

An important source of explicit positive examples is the
orthogonal specialization of the two-parameter biorthogonal
family introduced by Askey~\cite{A82}. Sri Ranga~\cite{SriRanga2010}
studied this specialization, obtaining Szeg\H{o} polynomials
associated with weights of the form
\[
    e^{\eta\varphi}\bigl[\sin(\varphi/2)\bigr]^{2s},
    \qquad 0<\varphi<2\pi,\quad
    s> -\tfrac12,\quad \eta\in\mathbb R,
\]
and established their hypergeometric representations and explicit
relations.
Vinet and Zhedanov~\cite[Sections~4.2--4.3]{VinetZhedanovAskey2021}
subsequently developed an algebraic interpretation of the
Askey biorthogonal polynomials through the meta-Jacobi algebra,
expressing their second-order differential equation as a
generalized eigenvalue problem for two fixed differential
operators. In the positive orthogonal specialization, this
equation is precisely of the form
\eqref{eq:intro-degree-dependent-equation} with $A(z)=z(1-z)$.

A different route to explicit OPUC was developed by
Cantero and Iserles~\cite{CanteroIserles2016}, who obtained
generating functions for polynomials with geometrically progressing
Schur parameters, including a generating function governed by a
functional differential equation.

Borrego--Morell and Sri Ranga~\cite{BorregoSriRanga2017} considered
second-order equations with degree-dependent polynomial coefficients,
allowing initially a leading coefficient $A_n$ of degree one or
two, with distinct roots in the quadratic case. They subsequently
studied the Gauss and confluent hypergeometric normal forms, using
complex affine changes of variable that could depend on the degree.
By means of raising and lowering operators,
they established conditional hypergeometric descriptions under
polynomial-solvability restrictions and an additional compatibility
hypothesis in Theorem~1.6. They conjectured that, apart from the
monomials, the circular hypergeometric family exhausts the positive
OPUC solutions in the Gauss case, and that only the monomials occur
in the normalized confluent case. These conjectures constitute a
direct antecedent of the classification problem considered here.

Further developments include the work of Garza, Garza, and
Huertas~\cite{GarzaGarzaHuertas2020}, who combined recurrence,
structure, and connection formulas to derive second-order equations
for OPUC associated with spectral transformations and coherent pairs.
Bracciali, Rampazzi, and Silva Ribeiro~\cite{BraccialiRampazziRibeiro2023}
classified semiclassical circle weights satisfying Pearson-type
equations involving polynomials of degree at most two and derived
structure relations and nonlinear equations for their Verblunsky
coefficients. They subsequently considered Pearson-type equations
with coefficients of degree at most three and obtained further
scalar differential equations, including equations for families
generalizing the circular Jacobi and modified Bessel
examples~\cite{BraccialiRampazziRibeiro2025}. Branquinho,
Foulqui\'e-Moreno, and Rampazzi~\cite{BranquinhoFoulquieRampazzi2026}
used Riemann--Hilbert methods to derive differential equations and
relations for recurrence coefficients of generalized semiclassical
OPUC. Such results develop differential structure from prescribed
properties of the weight. The degree bounds in a Pearson equation
for a weight are, however, distinct from the bounds imposed in
\eqref{eq:intro-degree-dependent-equation} on the equations for the
polynomials themselves.

A particularly valuable source on the theory of semiclassical OPUC
is provided by the unpublished lecture notes of Magnus~\cite{Magnus2013},
which contain an extensive treatment of generalized Jacobi polynomials
on the unit circle, including their differential equations, recurrence
relations, and asymptotic properties.

Another formulation of the circle analogue of Bochner's problem
arises from the CMV representation. Cantero, Moral, and
Vel\'azquez~\cite{CanteroMoralVelazquez2003} introduced the
five-diagonal matrices representing multiplication by $z$ in an
orthonormal Laurent polynomial basis; see also
Simon~\cite{SimonCMV2007}. These matrices play a role analogous to
that of Jacobi matrices in the real-line theory. In the spirit of
the bispectral problem of Duistermaat and
Gr\"unbaum~\cite{DuistermaatGrunbaum1986}, one asks whether the same
Laurent polynomials are eigenfunctions of a differential operator
in $z$. Gr\"unbaum and Vel\'azquez~\cite{GrunbaumVelazquez2017}
proved that the existence of a common nonscalar ordinary differential
operator of finite order diagonal on the CMV Laurent basis forces
the orthogonality measure to be normalized Lebesgue measure, even
if finitely many initial eigenfunction conditions are omitted.

Vinet and Zhedanov~\cite{VinetZhedanovCMV2026} obtained
nontrivial examples of CMV bispectrality for Jacobi polynomials
on the unit circle by allowing differential--reflection
operators. Our formulation concerns a different spectral
structure: ordinary monic polynomials satisfying scalar
second-order differential equations with degree-dependent
lower-order coefficients and a fixed leading coefficient,
without requiring a common differential operator diagonal
on a Laurent polynomial basis.

Our main result, Theorem~\ref{thm:circle-main-classification}, gives
an exhaustive classification for the four canonical leading
coefficients. In addition to the normalized Lebesgue measure
$d\varphi/(2\pi)$, the only possibilities are the circular
hypergeometric measures
\[
    d\nu_q(e^{i\varphi})
    =C_q e^{(\pi-\varphi)\operatorname{Im}q}
      |1-e^{i\varphi}|^{2\operatorname{Re}q}\frac{d\varphi}{2\pi},
    \qquad \operatorname{Re}q>-\tfrac12,\quad q\ne0,
\]
where $C_q>0$ is a normalizing constant, and the Bernstein--Szeg\H{o}
measures
\[
    d\mu_{m,a}(e^{i\varphi})
    =\frac{1-|a|^2}{|e^{im\varphi}-a|^2}\frac{d\varphi}{2\pi},
    \qquad m\geq1,\quad 0<|a|<1.
\]
The latter have exactly one nonzero Verblunsky coefficient,
$\alpha_{m-1}=\overline a$; for $m=1$ they are the Poisson measures.
More precisely, no positive measure of infinite support occurs
for $\theta=1$; for $\theta=z$ only the Lebesgue and Poisson
($m=1$) measures occur; for $\theta=z(1-z)$ the possible measures
are the Lebesgue, Poisson, and circular hypergeometric measures; and
for $\theta=z^2$ they are the Lebesgue measure and all the measures
$\mu_{m,a}$.  In particular, absolute continuity is a consequence
of the classification, not a hypothesis. The theorem also
identifies the corresponding monic polynomial sequences and the
admissible differential coefficients. Moreover, the classification
remains valid in all four canonical cases when the differential
equation is imposed only in sufficiently large degrees
(see Remark~\ref{rem:circle-eventual-classification}).

The conjecture of~\cite{BorregoSriRanga2017} for the normalized
confluent equation ($\theta=z$, $B(z,n)=c_n-z$) is thereby confirmed.
However, the earlier analysis of the general confluent equation
and the Gauss conjecture must both be completed by the Poisson
family. These measures satisfy the corresponding polynomial-solvability
restrictions, and their omission in both cases can be traced to
the exclusion of nontrivial finite Verblunsky tails with
$\alpha_0\ne0$ in Proposition~1.3 of that paper. Moreover,
the nontrivial Euler families for $\theta=z^2$ arise from resonance
of the differential equation. Thus the existence of a positive
orthogonal sequence must be distinguished from uniqueness of the
monic polynomial solution in each degree. The classification
does not require nonresonance, nonvanishing of the Verblunsky
coefficients, or the additional ladder-operator compatibility
hypothesis of the 2017 result.

In the present manuscript we follow a different approach: we start with the differential equations
and recover the orthogonality measures by using coefficient
comparison, the location of the polynomial zeros, and compatibility
with the Szeg\H{o} recurrence. Special attention is given to vanishing
Verblunsky coefficients and to possible degeneracies of the
polynomial solutions. The preliminary results establish the
polynomial-solvability criteria, the fixed-coefficient case, and
the relation with the 2017 conjectures; the subsequent analysis
treats the four canonical equations and yields the classification
summarized in Theorem~\ref{thm:circle-main-classification} and
Table~\ref{tab:circle-main-classification}.

The remainder of the paper is organized as follows.
Section~\ref{sec:preliminaries-results} collects the necessary
preliminaries and states the main classification theorem.
Section~\ref{sec:proof-results} is devoted to its proof, which is
organized according to the four canonical configurations of the
leading coefficient. The Hermite case is excluded in
Subsection~\ref{subsec:circle-constant-leading}, while the confluent,
Gauss, and double-zero cases are completely classified in
Subsections~\ref{Confluency}, \ref{Gauss}, and
\ref{sec:circle-double-zero}, respectively. Finally, Section~\ref{sec:circle-signed-measures} discusses the
Pearson equation approach and outlines the extension of the
classification to signed measures.

 \section{Preliminaries and statement of the results}
\label{sec:preliminaries-results}

Let $\mu$ be a probability measure with infinite support on
$\mathbb T=\{z:|z|=1\}$, and let $\Phi_n$ be its monic orthogonal
polynomials. Thus
\[
    \int_{\mathbb T}\Phi_n(z)z^{-k}\,d\mu(z)=0,
       \qquad 0\leq k<n.
\]
they satisfy the Szeg\H{o}  recurrence relation given by 
\begin{equation}
\label{eq:circle-Szego}
    \Phi_{n+1}=z\Phi_n+\delta_{n+1}\Phi_n^*,\qquad
    \delta_n=\Phi_n(0),\qquad
    \alpha_{n-1}=-\overline{\delta_n},\qquad |\delta_n|<1,
\end{equation}
where $\Phi_n^*(z)=z^n\overline{\Phi_n(1/\overline z)}$ and the parameters  $\alpha_n$, $n\geq0$, are  the
Verblunsky coefficients of $\mu$.
We use the standard facts that all zeros of $\Phi_n$ lie strictly
inside the unit disk and that the complete monic sequence determines
the probability measure; see~\cite{SimonOPUC2005}.

Borrego--Morell and Sri Ranga~\cite{BorregoSriRanga2017} studied
second-order equations with polynomial coefficients varying with
the degree and obtained conditional hypergeometric descriptions
using raising and lowering operators.  
Our problem is to determine the measure itself, using Darboux
polynomials and cofactors together with positivity and
\eqref{eq:circle-Szego}.

We consider the differential equations
\begin{equation}\label{pre}
 A(z)y''+B(z,n)y'+C_ny=0,\qquad n\geq1,
\end{equation}
where $A\not\equiv0$ is independent of $n$,
$\deg A\leq2$, $\deg_z B\leq1$, and $C_n\in\mathbb C$.
Our objective is to classify, up to a common invertible complex
affine change of variable, the positive probability measures
of infinite support on the unit circle whose monic orthogonal
polynomials satisfy \eqref{pre} in every degree.

As in Bochner's classification on the real line, affine changes
of variable are understood to act simultaneously on the
differential equation and the orthogonality measure. In the
present setting, such a transformation sends the unit circle
to another circle and preserves positivity, infinite support,
and polynomial orthogonality. Accordingly, the four canonical
configurations considered below are understood relative to
the circle supporting the orthogonality measure, and the
classification extends to all their simultaneous affine
equivalents.

Since $A$ is fixed, there are constants $\rho,\sigma,\kappa$,
with $\sigma\kappa\neq0$, independent of $n$, such that
\[
 \frac{A(\rho+\sigma w)}{\sigma^2}
       =\kappa\theta(w),
\]
where
\[
 \theta(w)=
 \begin{cases}
   w(1-w), & \deg A=2,\quad \operatorname{disc}(A)\neq0,\\
   w,      & \deg A=1,\\
   1,      & \deg A=0,\\
   w^2,    & \deg A=2,\quad \operatorname{disc}(A)=0.
 \end{cases}
\]
Thus the single affine substitution $z=\rho+\sigma w$,
followed by division by $\kappa$, transforms \eqref{pre} into
\begin{equation}\label{prepre}
 \theta(w)Y_n''(w)+(c_n-b_nw)Y_n'(w)
       +\lambda_nY_n(w)=0,
\end{equation}
where $Y_n(w)=y(\rho+\sigma w)$ and
\[
 c_n-b_nw=\frac{B_n(\rho+\sigma w)}{\kappa\sigma},
 \qquad \lambda_n=\frac{C_n}{\kappa}.
\]

If $\Phi_n$ are the original monic orthogonal polynomials,
their monic transforms
\[
 \Psi_n(w)=\sigma^{-n}\Phi_n(\rho+\sigma w)
\]
are orthogonal for the common positive probability measure
$\nu=(F^{-1})_{\#}\mu$, where $F(w)=\rho+\sigma w$.
This measure is supported on the fixed circle
$\Gamma=\{w:|\rho+\sigma w|=1\}$, and the transformed
polynomials satisfy
\begin{equation}\label{eq:circle-general-DE}
 \theta(w)\Psi_n''(w)+B(w,n)\Psi_n'(w)
       +\lambda_n\Psi_n(w)=0,
 \qquad B(w,n)=c_n-b_nw.
\end{equation}
The canonical equation is understood together with its
transformed orthogonality circle; the original measure is
recovered by $\mu=F_{\#}\nu$.

The dependence of $B_n$ on the degree is essential in the
linear and distinct-root quadratic classes. A fixed affine
coefficient $B$ gives no positive infinite-support
orthogonality measure on a circle, by the three-term
recurrence obstruction used in
Proposition~\ref{prop:circle-fixed-B}.

Leading-coefficient comparison gives the precise forms
\begin{align}
 \mathrm{(H)}\quad&
 \Phi_n''+(c_n-b_nz)\Phi_n'+nb_n\Phi_n=0.
 \label{eq:circle-H}\\
  \mathrm{(K)}\quad&
 z\Phi_n''+(c_n-b_nz)\Phi_n'+nb_n\Phi_n=0,
 \label{eq:circle-K}\\
 \mathrm{(G)}\quad&
 z(1-z)\Phi_n''+(c_n-b_nz)\Phi_n'
       +n(n-1+b_n)\Phi_n=0,
 \label{eq:circle-G}\\
 \mathrm{(E)}\quad&
 z^2\Phi_n''+(c_n-b_nz)\Phi_n'
       -n(n-1-b_n)\Phi_n=0,\qquad n\geq1.
\end{align}

We first determine the conditions under which the differential
equations above admit a monic polynomial solution of degree $n$.
The following proposition extends the criterion for polynomial
solvability in Proposition~1.1 of~\cite{BorregoSriRanga2017} to
all four cases and describes the solutions at the exceptional
parameter values.

\begin{proposition}[Polynomial solvability and resonance]
\label{prop:circle-polynomial-solvability}
Let $n\ge1$ be fixed, and let $b_n,c_n\in\mathbb C$.
Consider one of the equations $\mathrm{(K)}$, $\mathrm{(G)}$,
$\mathrm{(H)}$ and $\mathrm{(E)}$, with the spectral coefficient
as given above. A monic polynomial solution of degree $n$ exists
and is unique if and only if the corresponding condition holds:
\begin{align}
\mathrm{(K)},\ \mathrm{(H)}:\quad
    &b_n\ne0,\\
\mathrm{(G)}:\quad
    &b_n\notin\{-2n+2,\ldots,-n+1\},
    \label{eq:circle-nonresonance}\\
\mathrm{(E)}:\quad
    &b_n\notin\{n-1,\ldots,2n-2\}.
\end{align}
Under these nonresonance conditions, $c_n$ is arbitrary.

At the remaining values of $b_n$, which we call resonant, the
existence conditions and corresponding monic solutions are as
follows, with $\tau\in\mathbb C$ arbitrary.
\begin{enumerate}
\item[\textup{(K)}] If $b_n=0$, existence holds if and only if
$c_n=1-n$. All monic solutions of degree $n$ are
\[
    p(z)=z^n+\tau.
\]

\item[\textup{(G)}] If $b_n=1-n-r$ for some integer
$0\le r\le n-1$, existence holds if and only if
\[
    c_n\in\{-\ell:\ell\in\mathbb Z,\ r\le\ell\le n-1\}.
\]
In this case, all monic solutions of degree $n$ are
\[
    p(z)=P_n(z)+\tau Q_r(z),
\]
where $P_n$ is the monic solution of degree $n$ with
$[z^r]P_n=0$, and $Q_r$ is the unique monic polynomial
solution of degree $r$ of the same differential equation,
with $b_n,c_n,\lambda_n$ unchanged.

\item[\textup{(H)}] If $b_n=0$, existence holds if and only if
$n=1$ and $c_n=0$. The monic solutions are then $p(z)=z+\tau$.

\item[\textup{(E)}] If $b_n=n+r-1$ for some integer
$0\le r\le n-1$, existence holds if and only if $c_n=0$.
All monic solutions of degree $n$ are
\[
    p(z)=z^n+\tau z^r.
\]
\end{enumerate}
Thus, at resonant parameter values, the set of monic solutions
of degree $n$ is either empty or a one-dimensional complex
affine space.
\end{proposition}

Proposition~\ref{prop:circle-polynomial-solvability} concerns
each individual differential equation. In the classification
of measures below,  no uniqueness or nonresonance  condition on the monic solution is imposed. Existence is part of the hypothesis that the orthogonal polynomials satisfy the differential equation.

 The Bochner characterization on the unit circle was already
observed by Cachafeiro and Su\'arez
\cite[p.~405]{CachafeiroSuarez1997}, who noted that it
singles out the monomial sequence. For completeness, we state
the corresponding result with the explicit restrictions
on the differential coefficients.

\begin{proposition}[Fixed affine coefficients and the Euler exception]
\label{prop:circle-fixed-B}
Let $\mu$ be a positive probability measure of infinite support
on $\mathbb T$, and let $\Phi_n$ be its monic orthogonal
polynomials.  Suppose that \eqref{eq:circle-general-DE} holds
for every $n\geq1$, with
\[
 \theta\in\{1,z,z(1-z),z^2\},\qquad
 B(z,n)=c-bz,
\]
where $b,c\in\mathbb C$ are independent of $n$.
Then necessarily
\begin{equation}
\label{eq:circle-fixed-B-only-case}
 \theta=z^2,\qquad c=0,\qquad
 \Phi_n(z)=z^n,\qquad
 d\mu(e^{i\varphi})=\frac{d\varphi}{2\pi},
\end{equation}
and
\[
 \lambda_n=-n(n-1-b).
\]
Conversely, \eqref{eq:circle-fixed-B-only-case} satisfies
the differential equations for every fixed $b\in\mathbb C$
with these spectral coefficients.  In particular, there is
no such measure for $\theta=1$, $\theta=z$, or
$\theta=z(1-z)$.
\end{proposition}

 \subsection{Statement of the results}
\label{subsec:circle-main-classification}

We first fix the notation for the measures and polynomial families
appearing in the classification. We denote by
$\mu_{\mathrm{Leb}}$ the normalized Lebesgue measure on $\mathbb T$:
\begin{equation}
\label{eq:circle-main-Lebesgue-measure}
    d\mu_{\mathrm{Leb}}(e^{i\varphi})
      =\frac{d\varphi}{2\pi}.
\end{equation}
Its monic orthogonal polynomials are $\Phi_n(z)=z^n$,
and its Verblunsky coefficients are identically zero.
This is the free case of OPUC;
see~\cite[Example~1.6.1]{SimonOPUC2005}.

For an integer $m\ge1$ and $0<|a|<1$, we define
\begin{equation}
\label{eq:circle-main-Poisson-families}
    d\mu_{m,a}(e^{i\varphi})
       =\frac{1-|a|^2}{|e^{im\varphi}-a|^2}
          \frac{d\varphi}{2\pi}.
\end{equation}
These are Bernstein--Szeg\H{o} measures, obtained by taking
$p(z)=z^m-a$ in the classical construction. This class goes back to
Szeg\H{o}'s work~\cite{Szego1919}; see
also~\cite[Section~4]{SimonOneFoot2005} and
\cite{SimonOPUC2005}.
For $m=1$, we retain the notation $\mu_a=\mu_{1,a}$ for the
classical Poisson-kernel measure:
\begin{equation}
\label{eq:circle-main-Poisson-measure}
    d\mu_a(e^{i\varphi})
       =\frac{1-|a|^2}{|e^{i\varphi}-a|^2}
          \frac{d\varphi}{2\pi}.
\end{equation}
We denote the monic orthogonal polynomials of $\mu_{m,a}$ by
\begin{equation}
\label{eq:circle-main-binomial-polynomials}
    \Pi_n^{(m,a)}(z)=
    \begin{cases}
        z^n,&0\le n<m,\\
        z^{n-m}(z^m-a),&n\ge m.
    \end{cases}
\end{equation}

For $q\in\mathbb C\setminus\{0\}$ with
$\operatorname{Re}q>-\tfrac12$, we also define the circular
hypergeometric measure
\begin{equation}
\label{eq:circle-main-hypergeometric-family}
    d\nu_q(e^{i\varphi})
       =C_q e^{(\pi-\varphi)\operatorname{Im}q}
          |1-e^{i\varphi}|^{2\operatorname{Re}q}
          \frac{d\varphi}{2\pi},\qquad 0<\varphi<2\pi,
\end{equation}
where $C_q>0$ normalizes the measure. This family is the positive
orthogonality specialization of the hypergeometric biorthogonal
system introduced by Askey in his 1982 discussion~\cite{A82}
of Szeg\H{o}'s paper. Sri Ranga~\cite{SriRanga2010} studied its
positive orthogonality and explicit Szeg\H{o} polynomials,
and the associated second-order differential equations were
considered in~\cite{BorregoSriRanga2017}. The corresponding
monic orthogonal polynomials admit the regular representation
\begin{equation}
\label{eq:circle-main-hypergeometric-polynomials}
    P_n^{(q)}(z)=
       \frac{(q+\overline q+1)_n}{(q+1)_n}
       {}_2F_1(-n,q+1;q+\overline q+1;1-z).
\end{equation}
Here $(u)_n$ denotes the rising factorial, with $(u)_0=1$.

With this notation, the Verblunsky coefficients of the three
families are, for $j\ge0$,
\begin{equation}
\label{eq:circle-main-Verblunsky}
\begin{aligned}
    \alpha_j(\mu_{\mathrm{Leb}})&=0,\\
    \alpha_j(\mu_{m,a})&=
       \begin{cases}
           \overline a,&j=m-1,\\
           0,&j\ne m-1,
       \end{cases}\\
    \alpha_j(\nu_q)&=
       -\frac{(q)_{j+1}}{(\overline q+1)_{j+1}}.
\end{aligned}
\end{equation}
Thus $\mu_{m,a}$ has exactly one nonzero Verblunsky coefficient,
whereas every Verblunsky coefficient of $\nu_q$ is nonzero.
The two families are therefore disjoint for the stated
parameter ranges. Allowing $a=0$ in $\mu_{m,a}$ or $q=0$ in
$\nu_q$ recovers the normalized Lebesgue measure.

Our main theorem determines which of these families occur
in each of the four canonical differential equations, together
with their monic orthogonal polynomials and admissible
differential coefficients. Both resonant and nonresonant
equations are included.

\begin{theorem}[Classification in the canonical unit-circle configurations]
\label{thm:circle-main-classification}
Let $\mu$ be a positive probability measure of infinite support
on $\mathbb T$, and let $\Phi_n$ be its monic orthogonal
polynomials. Fix
\[
    \theta\in\{1,z,z(1-z),z^2\}.
\]
Suppose that, for every $n\ge1$, there exist coefficients
$b_n,c_n,\lambda_n\in\mathbb C$ such that
\begin{equation}
\label{eq:circle-main-DE}
    \theta(z)\Phi_n''(z)+(c_n-b_nz)\Phi_n'(z)
        +\lambda_n\Phi_n(z)=0.
\end{equation}
Then the possible pairs $(\theta,\mu)$ are precisely those
listed in Table~\ref{tab:circle-main-classification}, and the
corresponding polynomial formulas hold for every $n\ge1$.

No nonresonance or uniqueness assumption is imposed.

The spectral coefficient is necessarily
\begin{equation}
\label{eq:circle-main-spectral-coefficients}
    \lambda_n=
    \begin{cases}
        nb_n,&\theta=1\text{ or }\theta=z,\\
        n(n-1+b_n),&\theta=z(1-z),\\
        -n(n-1-b_n),&\theta=z^2.
    \end{cases}
\end{equation}
For every $n\ge2$, the last two columns of
Table~\ref{tab:circle-main-classification} give all
admissible choices of the differential coefficients.
At degree one, all admissible choices are
\begin{equation}
\label{eq:circle-main-degree-one}
    b_1\in\mathbb C,\qquad c_1=-b_1\Phi_1(0).
\end{equation}

In particular, all admissible measures are absolutely
continuous.
\end{theorem}

\begin{table}[H]
\centering
\small
\renewcommand{\arraystretch}{1.35}
\setlength{\tabcolsep}{4pt}
\begin{tabular}{ccccc}
\hline
$\theta(z)$ & Measure & $\Phi_n(z)$
    & $b_n\ (n\ge2)$ & $c_n\ (n\ge2)$\\
\hline
$1$ & none & -- & -- & --\\
\hline
$z$ & Lebesgue & $z^n$ & arbitrary & $1-n$\\
    & Poisson $\mu_a$ & $z^{n-1}(z-a)$
        & $n/a$ & $2-n$\\
\hline
$z(1-z)$ & Lebesgue & $z^n$ & arbitrary & $1-n$\\
         & Poisson $\mu_a$ & $z^{n-1}(z-a)$
             & $2-2n+n/a$ & $2-n$\\
         & $\nu_q$ & $P_n^{(q)}(z)$
             & $q+2-n$ & $1-n-\overline q$\\
\hline
$z^2$ & Lebesgue & $z^n$ & arbitrary & $0$\\
      & $\mu_{m,a}$ & $\Pi_n^{(m,a)}(z)$
          & $\begin{cases}
                \text{arbitrary},&n<m,\\
                2n-m-1,&n\ge m
             \end{cases}$ & $0$\\
\hline
\end{tabular}
\caption{Complete classification in the canonical
unit-circle configurations.
Here $m\ge1$, $0<|a|<1$, and $q\ne0$ with
$\operatorname{Re}q>-1/2$.
Poisson denotes $\mu_a=\mu_{1,a}$.
The polynomial formulas hold for every $n\ge1$,
whereas the last two columns describe all admissible
differential coefficients for $n\ge2$.
The admissible coefficients at degree one are given by
\eqref{eq:circle-main-degree-one}.}
\label{tab:circle-main-classification}
\end{table}

 \begin{remark}[Classification from sufficiently large degrees]
\label{rem:circle-eventual-classification}
For each canonical leading coefficient
$\theta\in\{1,z,z(1-z),z^2\}$, the classification in
Theorem~\ref{thm:circle-main-classification} remains valid under
the weaker hypothesis that there exist an integer $N_0\geq1$
and complex sequences $b_n,c_n,\lambda_n$ such that
\eqref{eq:circle-main-DE} holds for every $n\geq N_0$.

In particular, Table~\ref{tab:circle-main-classification}
still describes precisely all admissible measures and their
corresponding monic orthogonal polynomials in every degree,
including those in which no differential equation is initially
imposed. The restrictions on the differential coefficients apply
to the degrees in which the equation is imposed. Conversely,
every listed family admits a realization of
\eqref{eq:circle-main-DE} in all positive degrees.
\end{remark}

\subsection{The 2017 conjecture and its correction}

The consideration of differential equations with degree-dependent
polynomial coefficients has a natural precedent in the work of
Shohat~\cite{Shohat1939}. In the semiclassical setting, he showed
that the orthogonal polynomials associated with weights satisfying
a polynomial Pearson equation, together with suitable boundary
conditions, satisfy second-order differential equations whose
polynomial coefficients may depend on $n$, but whose degrees
are uniformly bounded independently of $n$.

In the unit-circle setting, Borrego--Morell and Sri
Ranga~\cite{BorregoSriRanga2017} considered equations of the form
\[
    A(z,n)y''+B(z,n)y'+C_ny=0,
\]
where $A$ is either quadratic with distinct zeros or linear,
$\deg_z B\leq1$, and $C_n\in\mathbb C$.
They subsequently studied the canonical Gauss and confluent
equations and obtained conditional hypergeometric descriptions
through raising and lowering operators.

In the present work, we retain the degree dependence of $B$
and $C_n$, but restrict the leading coefficient to a fixed
polynomial $A$. This restriction, which preserves the configuration
of the leading zeros, provides a natural framework for establishing
a unit-circle analogue of Bochner's classification.  In the spirit of Bochner's original characterization,
the independence of $A$ from $n$ identifies a common
principal differential structure for the entire polynomial
sequence, thereby specifying the notion of classicality
adopted in the present work.

  Following \cite[Th.~1.6]{BorregoSriRanga2017},  they conjectured that, in the Gauss case, under
\eqref{eq:circle-nonresonance}, the sequence of positive OPUC
would belong to the circular hypergeometric family
\[
    R_n(z;s,\beta)
       ={}_2F_1(-n,s+\beta+1;\beta-s+1-n;z),
    \qquad s>-\tfrac12,\quad \beta\in i\mathbb R,
\]
up to normalization and inclusion of the monomial sequence.
In our notation, this corresponds to the family $P_n^{(q)}$
with $q=s+\beta$, where
$s=\operatorname{Re}q>-\frac12$ and
$\beta=i\operatorname{Im}q\in i\mathbb R$.

The historical conjecture for the normalized confluent equation
is confirmed by Theorem~\ref{thm:classification}. However, the
earlier analysis is not exhaustive for the  confluent
equation~\eqref{eq:circle-K}, just as the Gauss conjecture is not
exhaustive. Indeed, both equations admit the Poisson measures
described in Theorems~\ref{thm:classification}
and~\ref{thm:circle-finite-G}, whose monic orthogonal polynomials are
\[
    \Phi_n(z)=z^{n-1}(z-a),\qquad 0<|a|<1.
\]
These polynomials satisfy \eqref{eq:circle-K} with
\[
    c_n=2-n,\qquad b_n=\frac{n}{a},
\]
and \eqref{eq:circle-G} with
\[
    c_n=2-n,\qquad b_n=2-2n+\frac{n}{a}.
\]
The latter coefficients satisfy \eqref{eq:circle-nonresonance}.
In both cases, the Verblunsky coefficients are
$\alpha_0=\overline a\ne0$ and $\alpha_j=0$ for $j\ge1$.

The omission of the Poisson family in both cases can be traced
to Proposition~1.3 of~\cite{BorregoSriRanga2017}, which excludes
nontrivial finite Verblunsky tails when $\alpha_0\ne0$.
This exclusion is not valid in general.
In the Gauss case, if $M$ denotes the degree of the remaining
core polynomial, the cancellation leading to (2.10) is valid
only for $1\le k\le M$. For $M=1$, the order-two comparison
used in (2.12) is unavailable, since $(-1)_2=0$.
It is precisely the case $M=1$ that gives rise to the Poisson
family in both the general confluent and Gauss equations.

 \begin{remark}
In~\cite{BorregoSriRanga2017}, equations (2.27) and (2.28)
should read
\[
    y=-\frac{v'}{v},\qquad v''+P_1v'+Q_1v=0.
\]
While $v=\varphi_1$ proves Lemma~2.1, Lemma~2.2 requires a
linearly independent solution. Writing
$\varphi_1(z)=\kappa_1(z-\zeta)$, with $0<|\zeta|<1$,
we have $c_1=\zeta$ in the normalized confluent case.
In the Gauss case, $c_1=b_1\zeta$, and Proposition~1.1
permits choosing $b_1=1$ whenever the degree-one equation
can be freely selected. The nonresonance conditions
\[
    1-c_1\notin\mathbb Z,\qquad
    1+c_1-b_1\notin\mathbb Z,
\]
the latter only in the Gauss case, are then satisfied,
guaranteeing an independent Frobenius solution with
$v'/v$ holomorphic in $\Omega_\theta$. Setting
\[
    B_1=-\frac{v'}{v},\qquad
    A_1=\frac{\varphi_1'v-\varphi_1v'}{v},
\]
Abel's identity ensures that $A_1$ is nonvanishing and that
$A_1'/A_1=B_1-P_1$. The remaining identities follow from
the Riccati equation, completing the proof of Lemma~2.2.
\end{remark}

\section{Proof of the results}
\label{sec:proof-results}

\begin{proof}(\emph{Of Proposition~\ref{prop:circle-polynomial-solvability})}
Write $\theta(z)=\theta_2z^2+\theta_1z+\theta_0$ and
$p(z)=\sum_{j=0}^n u_jz^j$, where $u_n=1$ and
$u_{n+1}=u_{n+2}=0$. Leading-coefficient comparison gives
\[
    \lambda_n=n\bigl(b_n-\theta_2(n-1)\bigr).
\]
Consequently, comparison of the remaining coefficients gives
\begin{equation}
\label{eq:circle-polynomial-solvability-recurrence}
    d_{j,n}u_j
    +(j+1)(c_n+\theta_1j)u_{j+1}
    +\theta_0(j+1)(j+2)u_{j+2}=0,
    \qquad 0\le j<n,
\end{equation}
where
\[
    d_{j,n}=(n-j)\bigl(b_n-\theta_2(n+j-1)\bigr).
\]
This is a triangular system for $u_{n-1},\ldots,u_0$.
It has a unique solution if and only if
$d_{j,n}\ne0$ for every $0\le j<n$. Indeed, a singular
square coefficient matrix gives either no solution or
infinitely many solutions. Substituting $\theta_2=0,-1,0,1$
in cases $\mathrm{(K)}$, $\mathrm{(G)}$, $\mathrm{(H)}$
and $\mathrm{(E)}$, respectively, gives the stated
nonresonance conditions.

For $\theta_0=0$, iteration of
\eqref{eq:circle-polynomial-solvability-recurrence} gives
the polynomial in the nonresonant case explicitly as
\[
    p(z)=\sum_{j=0}^{n}(-1)^{n-j}\binom{n}{j}
       \left(\prod_{k=j}^{n-1}
          \frac{c_n+\theta_1 k}{b_n-\theta_2(n+k-1)}
       \right)z^j,
\]
where an empty product equals one. In the nonresonant case
$\mathrm{(H)}$,
writing $h(x)=p(x+c_n/b_n)$ gives $h''-b_nxh'+nb_nh=0$.
Coefficient comparison in powers of $x$ then gives
\[
    p(z)=\sum_{k=0}^{\lfloor n/2\rfloor}
        \frac{(-1)^k n!}{2^k k!(n-2k)!\,b_n^k}
        \left(z-\frac{c_n}{b_n}\right)^{n-2k}.
\]
These finite expressions require no additional restriction
on $c_n$.

It remains to consider resonant values. In case $\mathrm{(K)}$,
if $b_n=0$, the coefficient of $z^{n-1}$ forces $c_n=1-n$.
The resulting equation is
\[
    zp''+(1-n)p'=0,
\]
so the coefficient of $z^{j-1}$ is $j(j-n)u_j$.
Hence $u_j=0$ for $1\le j<n$, and $u_0$ is arbitrary.
In case $\mathrm{(H)}$, $b_n=0$ forces $c_n=0$ by the
coefficient of $z^{n-1}$. The equation reduces to $p''=0$,
which admits a monic polynomial of degree $n\ge1$ precisely
when $n=1$.

In the resonant cases $\mathrm{(G)}$ and $\mathrm{(E)}$,
there is a unique index $r\in\{0,\ldots,n-1\}$ for which
$d_{r,n}=0$. Since $\theta_0=0$, the equations with
$j>r$ determine $u_{r+1},\ldots,u_{n-1}$ uniquely.
At $j=r$, the necessary and sufficient compatibility condition is
\[
    (r+1)(c_n+\theta_1r)u_{r+1}=0.
\]
Using the recurrence for the higher coefficients and $u_n=1$,
this condition is equivalent to
\[
    \prod_{k=r}^{n-1}(c_n+\theta_1k)=0.
\]
It gives the stated condition on $c_n$ in case $\mathrm{(G)}$,
where $\theta_1=1$, and gives $c_n=0$ in case $\mathrm{(E)}$,
where $\theta_1=0$. Once the compatibility condition holds,
$u_r$ is arbitrary and the remaining coefficients
$u_{r-1},\ldots,u_0$ are uniquely determined. In case
$\mathrm{(G)}$, this yields the asserted affine family, with
$P_n$ obtained by taking $u_r=0$ and $Q_r$ by subtracting
that solution from the one with $u_r=1$.
In case $\mathrm{(E)}$, all off-diagonal coefficients vanish
when $c_n=0$, giving $p=z^n+\tau z^r$.
\end{proof}

 \begin{proof}(\emph{Of Proposition~\ref{prop:circle-fixed-B})}
First, positive OPUC cannot satisfy an ordinary monic
three-term recurrence
\begin{equation}
\label{eq:circle-impossible-Favard}
 \Phi_{n+1}=(z-\beta_n)\Phi_n-\gamma_n\Phi_{n-1},
 \qquad \gamma_n\neq0,
\end{equation}
for all sufficiently large $n$.
Indeed, if $\delta_n=\Phi_n(0)=0$, the Szeg\H{o} recurrence
gives $\Phi_n=z\Phi_{n-1}$ and
$\Phi_n^*=\Phi_{n-1}^*$.  Comparing the two recurrences yields
\[
 \delta_{n+1}\Phi_{n-1}^*
       =(-\beta_nz-\gamma_n)\Phi_{n-1}.
\]
For $n\geq2$, evaluation at a zero of $\Phi_{n-1}$ in the
unit disk forces $\delta_{n+1}=0$, because
$\Phi_{n-1}^*$ has no zeros there.  The displayed identity
then forces $\gamma_n=0$, a contradiction.  Therefore
$\delta_n\neq0$ for every sufficiently large $n$.

For such indices, the Szeg\H{o} recurrence implies
\[
 \Phi_{n+1}
 =\left(z+\frac{\delta_{n+1}}{\delta_n}\right)\Phi_n
 -\frac{\delta_{n+1}}{\delta_n}
       (1-|\delta_n|^2)z\Phi_{n-1}.
\]
Consecutive polynomials are relatively prime when
$\delta_n\neq0$: a common zero would also be a zero of
$\Phi_{n-1}^*$, contradicting the disk location of OPUC
zeros.  Comparing with \eqref{eq:circle-impossible-Favard}
would therefore force $\Phi_n$ to divide the nonzero affine
polynomial
\[
 \gamma_n-\frac{\delta_{n+1}}{\delta_n}
                     (1-|\delta_n|^2)z,
\]
which is impossible for $n\geq2$.

We now examine the four leading coefficients.

\smallskip
\noindent\emph{The constant case $\theta=1$.}
Leading-coefficient comparison gives $\lambda_n=nb$.
If $b=0$, comparison at degree $n-1$ forces $c=0$,
and the remaining equation $\Phi_n''=0$ is impossible
for $n\geq2$.
If $b\neq0$, the monic solutions are shifted and scaled
Hermite polynomials, with ordinary recurrence coefficients
\[
 \beta_n=\frac cb,\qquad \gamma_n=\frac nb\neq0.
\]
The preceding obstruction excludes this case.

\smallskip
\noindent\emph{The linear case $\theta=z$.}
Here again $\lambda_n=nb$.
If $b=0$, the coefficient of $z^{n-1}$ requires
$n(n-1+c)=0$ for every $n$, which is impossible.
If $b\neq0$, the monic solutions are scaled Laguerre
polynomials, whose ordinary recurrence has
\[
 \gamma_n=\frac{n(n+c-1)}{b^2}.
\]
This coefficient is nonzero for every sufficiently large
$n$, giving the same contradiction.

\smallskip
\noindent\emph{The Gauss case $\theta=z(1-z)$.}
The spectral coefficient is $\lambda_n=n(n-1+b)$.
For all sufficiently large degrees the monic solution is
unique and is a monic shifted Jacobi polynomial with
parameters $c-1$ and $b-c-1$.  Its recurrence coefficient is
\[
 \gamma_n=
 \frac{n(n+c-1)(n+b-c-1)(n+b-2)}
 {(2n+b-2)^2(2n+b-1)(2n+b-3)}.
\]
It is nonzero eventually, which is again impossible.
The displayed classical recurrence identities hold by
algebraic continuation for complex parameters wherever
their denominators do not vanish.  With fixed parameters,
exceptional degrees are finite and do not affect the argument.

\smallskip
\noindent\emph{The double-zero case $\theta=z^2$.}
Now $\lambda_n=-n(n-1-b)$.  Write
\[
 \Phi_n(z)=\sum_{j=0}^n a_{n,j}z^j,
 \qquad a_{n,n}=1.
\]
Coefficient comparison gives
\begin{equation}
\label{eq:circle-fixed-B-Euler-coefficients}
 c(j+1)a_{n,j+1}
 +(j-n)(j+n-1-b)a_{n,j}=0,
 \qquad 0\leq j<n.
\end{equation}
For all sufficiently large $n$, every factor
$j+n-1-b$ in this range is nonzero.

Suppose first that $c\neq0$.  Put
$s_n=[z^{n-1}]\Phi_n$ and $\delta_n=\Phi_n(0)$.
Equation \eqref{eq:circle-fixed-B-Euler-coefficients} gives
\begin{equation}
\label{eq:circle-fixed-B-Euler-top-constant}
 s_n=\frac{nc}{2n-b-2},\qquad
 \delta_n=\frac{c^n}{\prod_{j=0}^{n-1}(n+j-1-b)}\neq0.
\end{equation}
Consequently, for large $n$,
\[
 |\delta_n|
 \leq\left(\frac{|c|}{n-1-|b|}\right)^n,
 \qquad
 n^2\delta_{n+1}\overline{\delta_n}\longrightarrow0.
\]
The coefficient of $z^n$ in the Szeg\H{o} recurrence gives
\[
 s_{n+1}-s_n=\delta_{n+1}\overline{\delta_n}.
\]
On the other hand, the first formula in
\eqref{eq:circle-fixed-B-Euler-top-constant} gives
\[
 s_{n+1}-s_n
 =-\frac{c(b+2)}{(2n-b)(2n-b-2)}.
\]
Multiplication by $n^2$ and passage to the limit force
$b=-2$.  In that case $s_n=c/2$ is constant, so the
Szeg\H{o} identity would give
$\delta_{n+1}\overline{\delta_n}=0$, contradicting
\eqref{eq:circle-fixed-B-Euler-top-constant}.
Thus $c\neq0$ is impossible.

If $c=0$, equation
\eqref{eq:circle-fixed-B-Euler-coefficients} implies
$a_{n,j}=0$ for $j<n$ in every sufficiently large degree.
Hence $\Phi_n=z^n$ eventually.  Inverse Szeg\H{o}
recurrence then gives $\Phi_n=z^n$ at every earlier degree
as well.  Orthogonality makes all nonzero Fourier
coefficients of $\mu$ vanish, so $\mu$ is the normalized Lebesgue measure.  Finally, direct substitution verifies the
converse for arbitrary fixed $b$.
\end{proof}

\subsection{Exclusion of the Hermite case}
\label{subsec:circle-constant-leading}

A nonzero constant leading coefficient can be normalized to one.
We therefore consider
\begin{equation}
\label{eq:circle-H-general}
    y''+(c_n-b_nz)y'+\lambda_n y=0,
\end{equation}
where $b_n,c_n,\lambda_n\in\mathbb C$ may depend on $n$.
We first identify its polynomial solutions and then show that
they cannot form the monic orthogonal sequence of a positive
measure of infinite support on the unit circle, even if the
equation is required only in sufficiently large degrees.

\begin{proposition}[Polynomial solutions in the constant case]
\label{prop:circle-H-polynomials}
Let $n\ge2$. Equation \eqref{eq:circle-H-general} has a monic
polynomial solution of degree $n$ if and only if
\[
    b_n\ne0,\qquad \lambda_n=nb_n.
\]
In this case the solution is unique. Set
\begin{equation}
\label{eq:circle-H-center-scale}
    m_n=\frac{c_n}{b_n},\qquad s_n^2=\frac1{b_n}.
\end{equation}
Then the solution is
\begin{equation}
\label{eq:circle-H-polynomial}
\begin{split}
    p_n(z)
      &=s_n^n\operatorname{He}_n\!\left(\frac{z-m_n}{s_n}\right)\\
      &=\sum_{j=0}^{\lfloor n/2\rfloor}
          \frac{(-1)^j n!\,s_n^{2j}}
               {2^j j!(n-2j)!}(z-m_n)^{n-2j},
\end{split}
\end{equation}
where $\operatorname{He}_n$ is the monic probabilists' Hermite
polynomial, given by $\operatorname{He}_n(x)=2^{-n/2}H_n(x/\sqrt{2})$, with $H_n$ as in \cite[Sect.~5.5, pp.~105--106]{Szego1975}.  The expression is independent of the choice of $s_n$.
\end{proposition}

\begin{proof}
Comparison of the coefficient of $z^n$ gives $\lambda_n=nb_n$.
If $b_n=0$, the coefficient of $z^{n-1}$ gives $c_n=0$,
leaving $p_n''=0$, which is impossible for $n\ge2$.
For $b_n\ne0$, put
$h(x)=s_n^{-n}p_n(m_n+s_nx)$. The equation becomes
\[
    h''-xh'+nh=0,
\]
whose monic degree-$n$ solution is $\operatorname{He}_n$.
Uniqueness also follows by subtracting two monic solutions:
if their nonzero difference has degree $j<n$, its highest
coefficient is multiplied by $(n-j)b_n\ne0$.
\end{proof}

At degree one, $p_1=z-m$ satisfies
\eqref{eq:circle-H-general} precisely when
$\lambda_1=b_1$ and $c_1=b_1m$.
The only solvable resonance in positive degree occurs at
$n=1$, with $b_1=c_1=\lambda_1=0$; in this case every monic
linear polynomial $p_1(z)=z+\tau$, $\tau\in\mathbb C$, is
a solution.
At degree zero the condition is
$\lambda_0=0$. These low-degree freedoms do not affect the
following exclusion result.

\begin{theorem}[No positive OPUC realization of the constant case]
\label{thm:circle-H-exclusion}
There is no positive probability measure $\mu$ of infinite
support on $\mathbb T$ whose monic orthogonal polynomials
$\Phi_n$ satisfy \eqref{eq:circle-H-general} for every
sufficiently large $n$, with arbitrary complex coefficients
$b_n,c_n,\lambda_n$.
Consequently, a nonzero constant leading coefficient contributes
no measures to the positive unit-circle classification.
No absolute continuity or additional nonresonance hypothesis
is required.
\end{theorem}

\begin{proof}
Suppose such a measure exists. Write
\[
    \delta_n=\Phi_n(0),\qquad
    \Phi_{n+1}=z\Phi_n+\delta_{n+1}\Phi_n^*,
    \qquad |\delta_n|<1\quad(n\ge1).
\]

Choose $N\ge4$ such that the differential equation holds for
$n\ge N$. Proposition~\ref{prop:circle-H-polynomials} gives
$b_n\ne0$, $\lambda_n=nb_n$, and the zeros of $\Phi_n$ are
\begin{equation}
\label{eq:circle-H-zeros}
    \xi_{j,n}=m_n+s_nx_{j,n},\qquad 1\le j\le n,
\end{equation}
where $x_{j,n}$ are the real zeros of $\operatorname{He}_n$.
By positivity of the measure, $|\xi_{j,n}|<1$.  Their symmetry and Newton's identities give
\begin{equation}
\label{eq:circle-H-real-moments}
\begin{gathered}
    \sum_jx_{j,n}=\sum_jx_{j,n}^3=0,\\
    \sum_jx_{j,n}^2=n(n-1),\qquad
    \sum_jx_{j,n}^4=n(n-1)(2n-3).
\end{gathered}
\end{equation}
For the even sums, use the coefficients
$[x^{n-2}]\operatorname{He}_n=-n(n-1)/2$ and
$[x^{n-4}]\operatorname{He}_n=n(n-1)(n-2)(n-3)/8$.

First we prove exponential decay of the constant terms.
Let $L_n=\max_j|x_{j,n}|$. By
\eqref{eq:circle-H-real-moments},
\[
    L_n^2\ge
       \frac{\sum_jx_{j,n}^4}{\sum_jx_{j,n}^2}
       =2n-3.
\]
Both $L_n$ and $-L_n$ are Hermite zeros. Since
$|m_n\pm s_nL_n|<1$, averaging their squared moduli gives
\[
    |m_n|^2+|s_n|^2L_n^2<1.
\]
Consequently,
\begin{align}
\label{eq:circle-H-mean-square}
    \frac1n\sum_j|\xi_{j,n}|^2
      &=|m_n|^2+(n-1)|s_n|^2\notag\\
      &\le |m_n|^2+
          \frac{n-1}{2n-3}(1-|m_n|^2)\notag\\
      &\le1-\frac{1-|m_n|^2}{3}.
\end{align}
The arithmetic--geometric mean inequality and $1-u\le e^{-u}$
therefore imply
\begin{equation}
\label{eq:circle-H-delta-bound}
    |\delta_n|=\prod_j|\xi_{j,n}|
      \le\exp\!\left[-\frac{n(1-|m_n|^2)}6\right].
\end{equation}

Let $a_{1,n}=[z^{n-1}]\Phi_n$. The Szeg\H{o} recurrence gives
\begin{equation}
\label{eq:circle-H-centroid-Szego}
    a_{1,n}=\delta_1+
        \sum_{j=1}^{n-1}\delta_{j+1}\overline{\delta_j}.
\end{equation}
For $n\ge N$, $a_{1,n}=-nm_n$. Hence
\[
    n|m_n|\le n-1+|\delta_1|,
    \qquad n(1-|m_n|^2)\ge1-|\delta_1|>0.
\]
It follows from \eqref{eq:circle-H-delta-bound} that
\[
    |\delta_n|\le q:=
      \exp\!\left[-\frac{1-|\delta_1|}{6}\right]<1
      \qquad(n\ge N).
\]
Substituting this uniform bound into
\eqref{eq:circle-H-centroid-Szego}, and absorbing its finitely
many initial terms into the error, yields
\[
    |m_n|\le q^2+O(n^{-1}).
\]
Choose $r$ with $q^2<r<1$. For all sufficiently large $n$,
$|m_n|\le r$, so another application of
\eqref{eq:circle-H-delta-bound} gives
\begin{equation}
\label{eq:circle-H-exponential-delta}
    |\delta_n|\le\eta^n,
    \qquad \eta=\exp[-(1-r^2)/6]\in(0,1).
\end{equation}

Next we show that the first four complex power sums converge
exponentially. For $1\le k\le4$ and $n\ge4$, put
$a_{k,n}=[z^{n-k}]\Phi_n$. Coefficient comparison in the
Szeg\H{o} recurrence gives
\[
    a_{k,n+1}-a_{k,n}
       =\delta_{n+1}\,
          \overline{[z^{k-1}]\Phi_n}.
\]
Since all zeros of $\Phi_n$ lie in $\mathbb D$,
\[
    |[z^{k-1}]\Phi_n|\le\binom{n}{k-1}.
\]
Thus, for any fixed $\tau\in(\eta,1)$, each $a_{k,n}$ has
a limit with error $O(\tau^n)$. Newton's identities give the
same conclusion for
\[
    T_{k,n}=\sum_j\xi_{j,n}^k,\qquad 1\le k\le4.
\]
Write
\begin{equation}
\label{eq:circle-H-UV}
    U_n=T_{1,n}=U+O(\tau^n),\qquad
    V_n=T_{2,n}=V+O(\tau^n).
\end{equation}
Equations \eqref{eq:circle-H-zeros} and
\eqref{eq:circle-H-real-moments} imply
\begin{equation}
\label{eq:circle-H-center-power-sums}
    m_n=\frac{U_n}{n},\qquad
    s_n^2=\frac{V_n-U_n^2/n}{n(n-1)}.
\end{equation}
The third and fourth power sums are consequently
\begin{align}
\label{eq:circle-H-third-sum}
    T_{3,n}
       &=\frac{3U_nV_n}{n}-\frac{2U_n^3}{n^2},\\
\label{eq:circle-H-fourth-sum}
    T_{4,n}
       &=\frac{U_n^4}{n^3}
         +\frac{6U_n^2}{n^2}
             \left(V_n-\frac{U_n^2}{n}\right)\notag\\
       &\qquad
         +\frac{2n-3}{n(n-1)}
             \left(V_n-\frac{U_n^2}{n}\right)^2.
\end{align}
Because $U_n,V_n$ are bounded, both expressions tend to zero.
Their exponential convergence therefore gives
\begin{equation}
\label{eq:circle-H-third-fourth-decay}
    T_{3,n}=O(\tau^n),\qquad T_{4,n}=O(\tau^n).
\end{equation}
Multiplying \eqref{eq:circle-H-fourth-sum} by $n$ and taking
limits yields $0=2V^2$, so $V=0$ and $V_n=O(\tau^n)$.
Multiplying \eqref{eq:circle-H-third-sum} by $n^2$ now yields
$0=-2U^3$. Thus $U=0$ and
\begin{equation}
\label{eq:circle-H-UV-decay}
    U_n=O(\tau^n),\qquad V_n=O(\tau^n).
\end{equation}

Finally, these estimates force the measure to be Lebesgue.
The reality of the $x_{j,n}$, their zero sum, and
$|s_n|^2=|s_n^2|$ give
\begin{align*}
    \sum_j|\xi_{j,n}|^2
       &=n|m_n|^2+n(n-1)|s_n|^2\\
       &=\frac{|U_n|^2}{n}
           +\left|V_n-\frac{U_n^2}{n}\right|
        =O(\tau^n).
\end{align*}
Hence
\[
    \sum_j|\xi_{j,n}|
       \le\sqrt{n\sum_j|\xi_{j,n}|^2}
       =O\!\left(\sqrt n\,\tau^{n/2}\right)\longrightarrow0.
\]
Using
\[
    \Phi_n^*(z)=\prod_{j=1}^n(1-\overline{\xi_{j,n}}z),
\]
we obtain
\[
    \sup_{|z|=1}|\Phi_n^*(z)-1|
       \le\exp\!\left(\sum_j|\xi_{j,n}|\right)-1
       \longrightarrow0.
\]
Conjugating the orthogonality relation for the test exponent
$n-k$ gives
\[
    \int_{\mathbb T}\Phi_n^*(z)z^{-k}\,d\mu(z)=0,
       \qquad 1\le k\le n.
\]
For each fixed $k\ge1$, passage to the limit shows that
$\int_{\mathbb T}z^{-k}\,d\mu=0$. Positivity gives the
conjugate moments, so uniqueness of Fourier coefficients of
finite measures implies
\[
    d\mu(e^{i\varphi})=\frac{d\varphi}{2\pi}.
\]
Its monic orthogonal polynomials are $\Phi_n=z^n$.
Substituting them into \eqref{eq:circle-H-general}, with
$\lambda_n=nb_n$, gives
\[
    n(n-1)z^{n-2}+nc_nz^{n-1}=0,
\]
which is impossible for $n\ge2$.
\end{proof}

The exclusion also applies to measures with atomic or
singular continuous components. By a common invertible
complex affine change of variable, the same conclusion
holds on any circle.

\subsection{Classification of the   confluent case $\theta(z)=z$}\label{Confluency}

We consider \eqref{eq:circle-K} with arbitrary $b_n\in\mathbb C$,
including the degenerate case $b_n=0$.

\subsubsection{A finite coefficient exclusion}

If $p(z)=\sum_{j=0}^n a_jz^j$, $a_n=1$, solves
\eqref{eq:circle-K} at degree $n$, its coefficients satisfy
\begin{equation}\label{eq:coefficient-recurrence}
 (j+1)(c_n+j)a_{j+1}+b_n(n-j)a_j=0,
 \qquad 0\leq j<n.
\end{equation}
When $b_n\ne0$ and $r=\operatorname{ord}_0 p$, the nonzero
coefficients are exactly those with $r\leq j\leq n$.
To see this, a zero coefficient $a_j$ with $r<j<n$ would force
$a_{j-1}=0$ by \eqref{eq:coefficient-recurrence}, and then inductively
$a_r=0$. If $r\geq1$, comparison at $j=r-1$ also gives
$c_n=1-r$.

\begin{lemma}[A nonzero constant term cannot persist at high degree]
\label{lem:finite-bound}
Suppose $\Phi_n$ and $\Phi_{n+1}$ satisfy \eqref{eq:circle-Szego} and
\eqref{eq:circle-K}, with $b_n b_{n+1}\ne0$. If $n\geq15$, then
$\delta_{n+1}=0$.
\end{lemma}

\begin{proof}
Suppose $\delta_{n+1}\ne0$. Abbreviate
\[
 b=b_n,\quad c=c_n,\quad B=b_{n+1},\quad C=c_{n+1},
 \qquad \Phi_n(z)=\sum_{j=0}^n a_jz^j.
\]
Every coefficient of $\Phi_{n+1}$ is nonzero by
\eqref{eq:coefficient-recurrence}; in particular,
\begin{equation}\label{eq:C-nonzero}
 C\ne0.
\end{equation}
Put $r=\operatorname{ord}_0\Phi_n$ and $m=n-r$. In the Szeg\H{o}
recurrence the term $z\Phi_n$ has nonzero coefficients precisely at
indices $r+1,\ldots,n+1$, while $\delta_{n+1}\Phi_n^*$ has them
precisely at indices $0,\ldots,m$. Since their sum has no zero
coefficient, necessarily $r\leq m$.

For the coefficient index $k$, define the rational functions
\begin{align*}
 R(k)&=-\frac{b(n+1-k)}{k(c+k-1)},\\
 S(k)&=-\frac{(n-k)(\bar c+n-k-1)}{\bar b(k+1)},\\
 T(k)&=-\frac{B(n+1-k)}{(k+1)(C+k)}.
\end{align*}
They are, respectively, the ratios of successive coefficients of
$z\Phi_n$, $\delta_{n+1}\Phi_n^*$, and $\Phi_{n+1}$ wherever the
corresponding coefficients are nonzero.

Introduce the two polynomials
\begin{align}
 A(k)&=Bk(c+k-1)-b(k+1)(C+k),\label{eq:A}\\
 D(k)&=(n-k)(\bar c+n-k-1)(C+k)
             -B\bar b(n+1-k).\label{eq:D}
\end{align}
The polynomial $D$ has degree three and leading coefficient one.
For $k=0,\ldots,r-1$, both coefficients being compared lie below
the support of $z\Phi_n$. Hence $S(k)=T(k)$, or equivalently
$D(k)=0$. Thus
\begin{equation}\label{eq:r-bound}
 r\leq3.
\end{equation}

Extend $a_j$ by zero outside $0\leq j\leq n$, and put
\[
 U_k=a_{k-1},\qquad Z_k=\delta_{n+1}\overline{a_{n-k}}.
\]
For $r+1\leq k\leq m$, the coefficients $U_k,Z_k,U_{k+1}$
are nonzero, and
\[
 U_{k+1}=R(k)U_k,\qquad Z_{k+1}=S(k)Z_k.
\]
The second identity also holds at $k=m$: here $Z_{m+1}=0$ and
\[
 S(m)=-\frac{r(\overline c+r-1)}{\overline b(m+1)}=0,
\]
using $c=1-r$ when $r\geq1$. The coefficient recurrence for
$\Phi_{n+1}$ therefore gives
\begin{equation}\label{eq:coefficient-sum}
 [R(k)-T(k)]U_k=[T(k)-S(k)]Z_k,
 \qquad r+1\leq k\leq m.
\end{equation}
Use this identity at $k$ and $k+1$. Eliminating $U_k,Z_k$ without
dividing by $A$ or $D$ gives
\[
 [T(k)-S(k)][R(k+1)-T(k+1)]R(k)
 -[R(k)-T(k)][T(k+1)-S(k+1)]S(k)=0.
\]
After clearing the nonzero denominators, this is
\begin{equation}\label{eq:H-zeros}
 H(k)=0,\qquad r+1\leq k\leq m-1,
\end{equation}
where
\begin{equation}\label{eq:H}
 \begin{split}
 H(k)={}&(c+k)(\bar c+n-k-1)D(k+1)A(k)\\
        &\hspace{8mm}-|b|^2D(k)A(k+1).
 \end{split}
\end{equation}
All denominators cleared here are nonzero in this index range:
this follows from the nonvanishing coefficient intervals and
\eqref{eq:coefficient-recurrence}. Thus no exceptional values of
$A(k)$ or $D(k)$ have been discarded.

By \eqref{eq:C-nonzero}, $A(0)=-bC\ne0$. Set $d=\deg A$, so
$0\leq d\leq2$. The first term in \eqref{eq:H} has degree $d+5$
and leading coefficient equal to minus the leading coefficient of
$A$; the second has degree $d+3$. Therefore $H$ is a nonzero
polynomial of degree at most seven. But \eqref{eq:H-zeros} supplies
\[
 m-r-1=n-2r-1\geq n-7\geq8
\]
distinct zeros, by \eqref{eq:r-bound} and $n\geq15$. This is a
contradiction.
\end{proof}

\begin{lemma}[Zero values of the scale parameter]
\label{lem:zero-b}
Assume that \eqref{eq:circle-K} holds for every degree
$k\geq N_0$, for some $N_0\geq1$. Under the positive-measure
assumptions, if $b_n=0$ for some
$n\geq\max\{3,N_0\}$, then the complete sequence consists
of monomials.
\end{lemma}

\begin{proof}
The equation at that degree gives
$c_n=1-n$ and $\Phi_n=z^n+\delta_n$. If $\delta_n\ne0$, put
$\varepsilon=\delta_{n+1}$. The Szeg\H{o} recurrence gives
\[
 \Phi_{n+1}=z^{n+1}+\delta_nz
                  +\varepsilon(1+\overline{\delta_n}z^n).
\]
If $b_{n+1}\ne0$ and $\varepsilon\ne0$, all coefficients must
be nonzero, contradicting the zero coefficient of $z^2$ for
$n\geq3$. If $b_{n+1}\ne0$ and $\varepsilon=0$, the coefficient
equation at $j=1$ similarly contradicts the nonzero coefficient
of $z$. If $b_{n+1}=0$, the polynomial must be
$z^{n+1}+\varepsilon$, again contradicting the coefficient of $z$.
Thus $\delta_n=0$ and $\Phi_n=z^n$.

The inverse Szeg\H{o} recurrence, using $|\delta_j|<1$, gives
monomials at all earlier degrees. If a first later nonzero
$\delta_k$ existed, then $\Phi_k=z^k+\delta_k$ with $k\geq4$.
Such a polynomial cannot solve the equation with $b_k\ne0$,
and the preceding argument rules out $b_k=0$. Hence no later
nonzero constant term exists either.
\end{proof}

\subsubsection{Complete classification}

\begin{theorem}
\label{thm:classification}
Let $\mu$ be a positive probability measure of infinite support on
$\mathbb T$. There exist an integer $N_0\geq1$ and complex sequences
$b_n,c_n$ such that its monic OPUC satisfy \eqref{eq:circle-K} for
every $n\geq N_0$ if and only if one of the following holds.

\begin{enumerate}
\item The measure is the normalized Lebesgue measure:
\[
    d\mu(e^{i\varphi})=\frac{d\varphi}{2\pi},
    \qquad \Phi_n(z)=z^n\quad(n\geq1).
\]
In every degree in which the differential equation is imposed,
the admissible parameters are $c_n=1-n$, with arbitrary $b_n$.

\item There exists a unique $a$ with $0<|a|<1$ such that
$\mu$ is the Poisson measure
\[
    d\mu(e^{i\varphi})=
    \frac{1-|a|^2}{|e^{i\varphi}-a|^2}\frac{d\varphi}{2\pi},
    \qquad
    \Phi_n(z)=z^{n-1}(z-a)\quad(n\geq1).
\]
In every imposed degree $n\geq2$, the admissible parameters
are necessarily
\[
    b_n=\frac{n}{a},\qquad c_n=2-n.
\]
If the equation is also imposed at degree one, the admissible
parameters are $c_1=ab_1$, with arbitrary $b_1$.
\end{enumerate}

In both cases, the polynomial formulas hold in every degree,
including those in which no differential equation is initially
imposed. Moreover, both families admit a realization of
\eqref{eq:circle-K} in all positive degrees. In particular, every admissible measure is absolutely
continuous with respect to Lebesgue measure. For the normalized confluent equation $b_n\equiv1$, even when
imposed only in sufficiently large degrees, the normalized
Lebesgue measure is the only possibility.
\end{theorem}

\begin{proof}
Assume that \eqref{eq:circle-K} holds for every $n\geq N_0$,
for some $N_0\geq1$. If $b_n=0$ for some
$n\geq\max\{3,N_0\}$, Lemma~\ref{lem:zero-b} gives the
monomial case. Otherwise, putting $N_*=\max\{15,N_0\}$,
we have $b_nb_{n+1}\ne0$ for every $n\geq N_*$.
Lemma~\ref{lem:finite-bound} therefore gives
\[
    \delta_k=0,\qquad k\geq N_*+1.
\]
In either case, only finitely many Verblunsky coefficients
can be nonzero.

If every $\delta_k=0$ for $k\geq1$, the Szeg\H{o} recurrence
yields the monomials. Substitution gives $c_n=1-n$, with
$b_n$ arbitrary, in every imposed degree.

Suppose instead that $M\geq1$ is the largest index with
$\delta_M\ne0$, and put $Q=\Phi_M$. Then
\[
    \Phi_n=z^{n-M}Q\quad(n\geq M),\qquad Q(0)\ne0.
\]
For every $n\geq\max\{N_0,M+1\}$, comparison of the
lowest nonzero coefficient in the differential equation gives
\[
    (n-M)(n-M-1+c_n)Q(0)=0.
\]
Consequently,
\[
    c_n=M+1-n.
\]
Substituting this value and canceling $z^{n-M}$ gives
\begin{equation}
\label{eq:tail-Q}
    zQ''+(n-M+1-b_nz)Q'+Mb_nQ=0.
\end{equation}
Subtracting the equations at $n+1$ and $n$ yields
\[
    Q'+(b_{n+1}-b_n)(MQ-zQ')=0.
\]
Since $Q$ has positive degree, $\Delta=b_{n+1}-b_n\ne0$.
Consequently, for $a=1/\Delta\ne0$,
\[
    (z-a)Q'=MQ.
\]
The monic solution is $Q=(z-a)^M$.

Although the original differential equation need not be imposed
at degree $M$, the reduced equation \eqref{eq:tail-Q}
is sufficient to exclude $M\geq2$.
Indeed, if $M\geq2$, substituting $Q=(z-a)^M$ into
\eqref{eq:tail-Q} and canceling $(z-a)^{M-2}$ gives
\[
    M(M-1)z
    +M(n-M+1-b_nz)(z-a)
    +Mb_n(z-a)^2=0.
\]
Evaluating at $z=a$ yields
\[
    aM(M-1)=0,
\]
contrary to $a\ne0$ and $M\geq2$. Hence $M=1$.

Thus
\[
    \Phi_n(z)=z^{n-1}(z-a),\qquad n\geq1,
\]
including all degrees in which no differential equation was
initially imposed. Since $\delta_1=-a$, positivity gives
$0<|a|<1$.

In every imposed degree $n\geq2$, the preceding reduction
gives $c_n=2-n$. Substituting $Q=z-a$ and $M=1$ into
\eqref{eq:tail-Q} gives $n-ab_n=0$, and consequently
\[
    b_n=\frac{n}{a},\qquad c_n=2-n.
\]
If the equation is imposed at degree one, substitution of
$\Phi_1=z-a$ gives $c_1=ab_1$, with $b_1$ arbitrary.

Finally, the normalized Lebesgue measure has the monomials
as its OPUC. For the displayed Poisson measure,
\[
    \int_{\mathbb T}z^j\,d\mu(z)=a^j,\qquad j\geq0.
\]
Hence, for $0\leq k<n$,
\[
    \int_{\mathbb T}z^{n-1}(z-a)z^{-k}\,d\mu(z)
    =a^{n-k}-a\,a^{n-1-k}=0.
\]
These are therefore its monic OPUC. The complete monic sequence
determines the probability measure, which identifies the measure
in both cases.

Conversely, direct substitution verifies that both families
satisfy \eqref{eq:circle-K} in every positive degree with the
admissible parameter choices stated in the theorem.
Their explicit densities establish the asserted absolute
continuity.

Finally, if $b_n=1$ for every sufficiently large $n$, the
Poisson case would require $a=n$ in every such degree,
contrary to $|a|<1$. Thus only the normalized Lebesgue
measure occurs in the normalized confluent case.
\end{proof}

\subsection{Classification of the Gauss case $\theta(z)=z(1-z)$}\label{Gauss}

We consider \eqref{eq:circle-G} with arbitrary complex
parameters, including the solvable resonant cases described in
Proposition~\ref{prop:circle-polynomial-solvability}.
It suffices to assume  that \eqref{eq:circle-G} holds for every
$n\ge N_0$, for some integer $N_0\ge1$. Set
\[
    N_1=\max\{6,N_0+1\},\qquad
    N_2=\max\{10,N_0+1\}.
\]
We distinguish two cases according to whether the set
\[
    \{j\ge0:\alpha_j\ne0\}
\]
is finite or infinite. The finite case is treated directly,
without a nonresonance hypothesis. In the infinite case,
we first prove that nonresonance holds for every degree
$n\ge N_1$. We then apply the nonresonant rigidity argument
at degrees $n\ge N_2$ and use the inverse Szeg\H{o} recurrence
to determine the complete sequence.

\subsubsection{The finite Verblunsky case}

The following argument applies to all parameters for which
the prescribed orthogonal polynomials satisfy \eqref{eq:circle-G}
for $n\ge N_0$, including resonant ones.

\begin{theorem}
\label{thm:circle-finite-G}
Assume \eqref{eq:circle-G} for every $n\ge N_0$ and that
only finitely many Verblunsky coefficients are nonzero.
Then $\mu$ is a Poisson measure, including Lebesgue as $a=0$:
\begin{equation}
\label{eq:circle-Poisson}
    d\mu_a(e^{i\varphi})
       =\frac{1-|a|^2}{|e^{i\varphi}-a|^2}
                       \frac{d\varphi}{2\pi},
       \qquad |a|<1.
\end{equation}
For $a\neq0$ its monic polynomials and one admissible
choice of differential coefficients in each degree $n\ge1$ are
\begin{equation}
\label{eq:circle-Poisson-DE}
    \Phi_n=z^{n-1}(z-a),\qquad
    c_n=2-n,\qquad b_n=2-2n+\frac n a.
\end{equation}
The displayed coefficients satisfy \eqref{eq:circle-nonresonance}
and are the only admissible coefficients for $n\ge2$.
In particular, they coincide with the prescribed coefficients
for $n\ge\max\{2,N_0\}$. At degree one, $b_1\in\mathbb C$
may be chosen arbitrarily and $c_1=ab_1$.
For Lebesgue, the admissible coefficients are $c_n=1-n$ and
arbitrary $b_n$ in every degree.
Conversely, every measure \eqref{eq:circle-Poisson} admits
\eqref{eq:circle-G} for all $n\ge1$.
\end{theorem}

\begin{proof}
The finite-tail hypothesis and \eqref{eq:circle-Szego}
give either the monomials or
$\Phi_{m+k}=z^kQ$ for all $k\geq0$, where
$m\geq1$, $\deg Q=m$, and $Q(0)\neq0$.
In the monomial case, the measure is the normalized Lebesgue measure, and the substitution
$\Phi_n=z^n$ in \eqref{eq:circle-G} gives $c_n=1-n$, with
$b_n$ arbitrary, in any degree in which the equation is imposed.
These coefficients also provide an extension to all earlier
degrees.
For the remaining case and $k\geq\max\{1,N_0-m\}$, substitution in
\eqref{eq:circle-G}
first gives $c_{m+k}=1-k$ from the lowest power at zero.
The resulting equation for $Q$ is
\[
    z(1-z)Q''+(k+1-\eta_kz)Q'
           +m(m-1+\eta_k)Q=0,\qquad
    \eta_k=b_{m+k}+2k.
\]
Subtract the equations for $k+1$ and $k$:
\[
    (1-\Delta_kz)Q'+m\Delta_kQ=0,\qquad
    \Delta_k=\eta_{k+1}-\eta_k.
\]
The case $\Delta_k=0$ would give $Q'=0$.
Otherwise this first-order equation forces
$Q=(z-a)^m$ for $a=1/\Delta_k$.
Here $0<|a|<1$. Such a point is ordinary for the
second-order equation for $Q$. A nonzero solution cannot
have both $Q(a)=Q'(a)=0$, so $m=1$.
It follows that $\Phi_n=z^{n-1}(z-a)$ for every $n\ge1$.
For this recovered polynomial, a Gauss equation in any degree
$n\ge2$ requires, by coefficient comparison at $z^{n-2}$
and $z^{n-1}$,
\[
    c_n=2-n,\qquad a(2n-2+b_n)=n.
\]
These conditions are also sufficient, giving precisely
the coefficients in
\eqref{eq:circle-Poisson-DE}; they are necessary in every
prescribed degree $n\ge2$ and admissible also in the earlier
degrees. At degree one, substitution of $\Phi_1=z-a$
gives $c_1=ab_1$, which includes the displayed choice.

Orthogonality to $1$ gives
$\int z^n\,d\mu=a\int z^{n-1}\,d\mu$, hence the moments
are $a^n$. Together with their conjugates these determine
\eqref{eq:circle-Poisson}.
Conversely, the Poisson moments are $a^n$, and
$\int\Phi_nz^{-j}\,d\mu_a=0$ for $0\leq j<n$.
Substitution verifies the differential equation.
Finally, a forbidden value of $b_n$ would require $n/a$
to be one of $0,1,\ldots,n-1$, whereas $|n/a|>n$.
\end{proof}

\subsubsection{The infinite Verblunsky case}
\label{subsec:circle-infinite-G}

Suppose now that infinitely many Verblunsky coefficients
are nonzero. In particular, the measure is different from
Lebesgue measure. No absolute-continuity assumption is imposed.

\paragraph{Reduction of the resonant equations.}
\label{subsec:circle-G-resonant}

We first establish auxiliary results valid for arbitrary
Verblunsky sequences, independently of the number of their
nonzero coefficients. They will ensure the nonresonance
required by the subsequent argument in the infinite case.

For each degree at which the differential equation is imposed,
we use the notation
\begin{equation}
\label{eq:circle-G-shifted-parameters}
    t_n=n-1+b_n,\qquad d_n=c_n+n-1,\qquad
    \delta_n=\Phi_n(0).
\end{equation}
The notation $\delta_n=\Phi_n(0)$ is used in all degrees,
including those below $N_0$.
The nonresonance condition \eqref{eq:circle-nonresonance}
is equivalent to
\begin{equation}
\label{eq:circle-G-t-nonresonance}
    t_n\notin\{0,-1,\ldots,1-n\}.
\end{equation}
For $\Phi_n(z)=\sum_{k=0}^n a_kz^k$, $a_n=1$,
equation \eqref{eq:circle-G} gives
\begin{equation}
\label{eq:circle-G-coefficient-recurrence}
    (k+1)(c_n+k)a_{k+1}
       +(n-k)(t_n+k)a_k=0,\qquad 0\le k<n.
\end{equation}

\begin{lemma}[Coefficient structure]
\label{lem:circle-G-coefficient-structure}
Let $\Phi_n$ be a monic polynomial of degree $n$ with all its
zeros in $\mathbb D$, satisfying \eqref{eq:circle-G}.
If the equation is nonresonant, its nonzero coefficients occur
precisely in degrees $\ell,\ell+1,\ldots,n$, where
$\ell=\operatorname{ord}_0\Phi_n$ denotes the order of
vanishing of $\Phi_n$ at zero.
If $t_n=-r$ for an integer $0\le r<n$, then
\begin{equation}
\label{eq:circle-G-resonant-form}
    c_n=1-n,\qquad
    \Phi_n(z)=z^n+\delta_n\sum_{k=0}^r
        (-1)^k\frac{n}{n-k}\binom rk z^k.
\end{equation}
Consequently, in either case:
\begin{enumerate}
\item if $\delta_n=0$, the nonzero coefficients occur in
consecutive degrees ending at $n$;
\item if $\delta_n\ne0$ and $a_{n-1}\ne0$, every coefficient
of $\Phi_n$ is nonzero.
\end{enumerate}
\end{lemma}

\begin{proof}
In the nonresonant case, \eqref{eq:circle-G-coefficient-recurrence}
shows that a vanishing coefficient forces all coefficients of
lower degree to vanish.

Suppose now that $t_n=-r$, and write
$Y(x)=\Phi_n(1-x)=\sum_{k=0}^n v_kx^k$. This polynomial satisfies
\[
    x(1-x)Y''+(b_n-c_n-b_nx)Y'+nt_nY=0.
\]
By the Gauss--Lucas theorem, $\Phi_n^{(k)}(1)\ne0$ for
$0\le k\le n$, since all zeros of $\Phi_n$ and of its
nonconstant derivatives lie strictly inside the unit disk.
Thus $v_k=(-1)^k\Phi_n^{(k)}(1)/k!\ne0$.
Coefficient comparison at $x^r$ gives
\[
    (r+1)(b_n-c_n+r)v_{r+1}=0.
\]
Hence $c_n=b_n+r=1-n$. Iterating
\eqref{eq:circle-G-coefficient-recurrence} from $a_0=\delta_n$
then gives
\[
    a_k=(-1)^k\frac{n}{n-k}\binom rk\delta_n,
    \qquad 0\le k<n,
\]
where $\binom rk=0$ for $k>r$. This proves
\eqref{eq:circle-G-resonant-form}.
In particular, a resonant polynomial with $\delta_n=0$ is
$z^n$, while a nonzero constant term and a nonzero coefficient
of $z^{n-1}$ force $r=n-1$ and all coefficients to be nonzero.
\end{proof}

\begin{lemma}[Exclusion of binomials and automatic nonresonance]
\label{lem:circle-G-automatic-nonresonance}
Suppose that the monic OPUC satisfy \eqref{eq:circle-G}
for every $n\ge N_0$. Then the following assertions hold:
\begin{enumerate}
\item For every $m\ge\max\{3,N_0\}$, a polynomial of
the form
\[
    \Phi_m(z)=z^m+\delta_m,\qquad \delta_m\ne0,
\]
cannot occur.

\item If $\Phi_m(z)=z^m$ for some
$m\ge\max\{2,N_0\}$, then the measure is the normalized Lebesgue  measure.

\item If the measure is not Lebesgue measure, then
\eqref{eq:circle-G-t-nonresonance} holds for every
$n\ge N_1$, where
\[
    N_1=\max\{6,N_0+1\}.
\]
\end{enumerate}
\end{lemma}

\begin{proof}
First, a polynomial $\Phi_m=z^m+\delta_m$ with
$m\ge\max\{3,N_0\}$ and $\delta_m\ne0$ cannot occur
in such a sequence.
Indeed, the Szeg\H{o} recurrence gives
\[
    \Phi_{m+1}=z^{m+1}+\delta_m z
           +\delta_{m+1}(1+\overline{\delta_m}z^m).
\]
If $\delta_{m+1}=0$, its nonzero coefficients at degrees
$1$ and $m+1$ are separated by vanishing coefficients.
If $\delta_{m+1}\ne0$, its constant and penultimate
coefficients are nonzero, yet the coefficient of $z^2$
vanishes. Since $m+1\ge N_0$, both alternatives contradict
Lemma~\ref{lem:circle-G-coefficient-structure} at degree $m+1$.

It follows that a monomial $\Phi_m=z^m$ at any degree
$m\ge\max\{2,N_0\}$ forces all subsequent constant terms
to vanish: otherwise the first subsequent nonzero constant
term would give the forbidden binomial at a degree
at least $\max\{3,N_0\}$. The inverse Szeg\H{o}
recurrence then gives $\Phi_j=z^j$ at all earlier degrees too,
and the measure is the normalized Lebesgue measure.

Now let $n\ge N_1$, suppose that the equation at degree $n$ is
resonant, and assume that the measure is not the normalized Lebesgue measure.
Since $n\ge6$ and $n-1\ge N_0$, the differential equations
at degrees $n$ and $n-1$ are both available.
By \eqref{eq:circle-G-resonant-form} and the preceding argument,
we have
\[
    \delta:=\delta_n\ne0,\qquad 1\le r\le n-1.
\]
Put $s=n-r-1$, $\gamma=|\delta|^2$, and
\[
    h_k=(-1)^k\frac{n}{n-k}\binom rk,\qquad 0\le k<n.
\]
The inverse Szeg\H{o} recurrence reads
\[
    z\Phi_{n-1}
       =\frac{\Phi_n-\delta\Phi_n^*}{1-\gamma}.
\]
Writing $\Phi_{n-1}(z)=\sum_{j=0}^{n-1}u_jz^j$, we obtain
\begin{equation}
\label{eq:circle-G-resonant-inverse}
    u_{k-1}=\frac{\delta h_k-\gamma h_{n-k}}{1-\gamma},
    \qquad 1\le k\le n-1.
\end{equation}
If $r<s$, then $h_1\ne0$ and
$h_{n-1}=h_{r+1}=h_s=0$, so
\[
    u_0\ne0,\qquad u_{n-2}\ne0,\qquad u_r=0.
\]
This contradicts Lemma~\ref{lem:circle-G-coefficient-structure}
at degree $n-1$. Therefore
\begin{equation}
\label{eq:circle-G-resonant-s-bound}
    r\ge s,\qquad s\le\frac{n-1}{2}.
\end{equation}

Using the identity
\[
    h_k=(-1)^k\binom nk
          \frac{\binom{n-1-k}{s}}{\binom{n-1}{s}},
\]
we may rewrite \eqref{eq:circle-G-resonant-inverse} as
\[
    u_{k-1}
       =\frac{(-1)^k\binom nk}
              {(1-\gamma)\binom{n-1}{s}}F(k),
       \qquad 1\le k\le n-1,
\]
where
\[
    F(x)=\delta\binom{n-1-x}{s}
             -(-1)^n\gamma\binom{x-1}{s}.
\]
Here the binomial coefficients with variable upper argument
are understood as polynomials, so $\deg F\le s$.

Set $C=c_{n-1}$ and $T=t_{n-1}$. The coefficient recurrence
\eqref{eq:circle-G-coefficient-recurrence} at degree $n-1$ gives
\[
    k(C+k-1)F(k+1)
       =(k+1)(T+k-1)F(k),\qquad 1\le k\le n-2.
\]
Consequently the polynomial
\[
    H(x)=x(C+x-1)F(x+1)
           -(x+1)(T+x-1)F(x)
\]
vanishes at $1,\ldots,n-2$. Its terms of highest degree
cancel, giving $\deg H\le s+1$. By
\eqref{eq:circle-G-resonant-s-bound} and $n\ge6$,
\[
    s+1\le\frac{n+1}{2}<n-2.
\]
Hence $H\equiv0$. Moreover,
\[
    F(0)=\delta\binom{n-1}{s}
             -(-1)^{n+s}\gamma\ne0,
\]
because $|\delta|\binom{n-1}{s}\ge|\delta|>|\delta|^2$.
Evaluating $H$ at $0$ and $-1$ yields $T=1$ and $C=2$,
respectively. The monic coefficient recurrence at degree
$n-1$ would then imply
\[
    \Phi_{n-1}(0)=(-1)^{n-1}
        \prod_{j=0}^{n-2}\frac{j+2}{j+1}
       =(-1)^{n-1}n,
\]
contrary to $|\Phi_{n-1}(0)|<1$. This excludes resonance
at every $n\ge N_1$ for a non-Lebesgue measure.
\end{proof}

\begin{lemma}[Consecutive zero constant terms]
\label{lem:circle-G-zero-tail}
Assume \eqref{eq:circle-G} for every $n\ge N_0$. If
$\delta_n=\delta_{n+1}=0$ for some $n\ge\max\{2,N_0\}$,
only finitely many Verblunsky coefficients are nonzero.
\end{lemma}

\begin{proof}
In view of Lemma~\ref{lem:circle-G-automatic-nonresonance},
we may assume that $\Phi_n$ is not a monomial.
Put $r=\operatorname{ord}_0\Phi_n\ge1$. Since
$\Phi_{n+1}=z\Phi_n$, the lowest nonzero powers in their
differential equations give $c_n=1-r$ and $c_{n+1}=-r$.
Subtracting the equations after substitution gives
\begin{equation}
\label{eq:circle-G-zero-first-order}
    z(1-Ez)\Phi_n'+(c_{n+1}+nEz)\Phi_n=0,
    \qquad E=2+b_{n+1}-b_n.
\end{equation}
Every nonzero zero of $\Phi_n$ is simple, since it is an
ordinary point of its second-order differential equation.
Equation \eqref{eq:circle-G-zero-first-order} therefore permits
at most one such zero. Hence
\[
    \Phi_n=z^{n-1}(z-a),\qquad 0<|a|<1.
\]
The zero constant term at degree $n+1$ extends this form
one degree further. If a subsequent constant term became
nonzero for the first time at degree $j\ge n+2$, then
\[
    \Phi_j=z^{j-1}(z-a)+\delta_j(1-\overline a z).
\]
Since $j\ge4$, its constant and penultimate coefficients
are nonzero, yet its coefficient of $z^2$ vanishes.
As $j\ge N_0$, this contradicts
Lemma~\ref{lem:circle-G-coefficient-structure} at degree $j$.
Hence all subsequent constant terms vanish.
\end{proof}

Returning to the infinite Verblunsky case, we conclude from
Lemma~\ref{lem:circle-G-automatic-nonresonance} that
\eqref{eq:circle-G-t-nonresonance} holds for every $n\ge N_1$.
Moreover, Lemma~\ref{lem:circle-G-zero-tail} excludes
consecutive zero constant terms $\delta_n=\delta_{n+1}=0$
for $n\ge\max\{2,N_0\}$.

\paragraph{The nonresonant rigidity argument.}

We may now apply the following coefficient-compatibility
and rigidity lemmas to the consecutive degrees $n,n+1$ with
$n\ge N_2$, since $n\ge10$ and both degrees lie in the
nonresonant range established above.
No differential equation is required at degrees below $N_0$,
and no nonresonance assumption on the original sequence is needed.
Once the polynomials at degrees $n\ge N_2$ have been
identified, the inverse Szeg\H{o} recurrence determines
all earlier degrees as well.

\begin{lemma}[A rational shift criterion]
\label{lem:circle-G-rational-shift}
Let $n\geq1$ be an integer, let $c,t\in\mathbb C$, and let
$g\in\mathbb C(x)\setminus\{0\}$ have a numerator and a denominator
of degree at most $h$.  Suppose that
\begin{equation}\label{eq:circle-G-rational-shift}
 \frac{g(x+1)}{g(x)}=
 \frac{(x+t-1)(\overline t+n-x-1)}
      {(x+c-1)(\overline c+n-x-1)}.
\end{equation}
Then at least one of the following holds:
\begin{equation}\label{eq:circle-G-small-shifts}
 t-c\in\mathbb Z,\quad |t-c|\leq h,
 \qquad\hbox{or}\qquad
 t+\overline c+n-2\in\mathbb Z,\quad
 |t+\overline c+n-2|\leq h.
\end{equation}
\end{lemma}

\begin{proof}
The zeros of the right-hand side of
\eqref{eq:circle-G-rational-shift} are
$a=1-t$ and $n+\overline t-1$, and its poles are
$p_1=1-c$ and $p_2=n+\overline c-1$, with multiplicities.
If $a$ is also a pole, one of the alternatives in
\eqref{eq:circle-G-small-shifts} holds with value zero.

Otherwise, write $\nu(\zeta)=\operatorname{ord}_{\zeta}g$.
The order of $g(x+1)/g(x)$ at $a$ is
$\nu(a+1)-\nu(a)>0$.
If $\nu(a)\geq0$, then $\nu(a+1)>0$.  Moving forward along
$a+\mathbb Z$, let $a+L$ be the first pole of the right-hand
side encountered.  Such a pole exists because $\nu$ has finite
support.  At $a+1,\ldots,a+L$ the orders of $g$ are positive,
so $1\leq L\leq h$.
If $\nu(a)<0$, move backward to the first pole $a-L$ of the
right-hand side.  At $a,a-1,\ldots,a-L+1$ the orders of $g$ are
negative, and again $1\leq L\leq h$.
The pole reached is $p_1$ or $p_2$.  Its difference from $a$
is respectively $t-c$ or $t+\overline c+n-2$, proving the claim.
\end{proof}

\begin{lemma}[Coefficient compatibility]
\label{lem:circle-G-compatibility}
Let $n\geq10$.  Suppose that two consecutive monic OPUC
$\Phi_n,\Phi_{n+1}$ satisfy \eqref{eq:circle-G}, nonresonance,
and \eqref{eq:circle-Szego}.  Assume that
\[
 \delta_{n+1}\neq0,\qquad
 r:=\operatorname{ord}_0\Phi_n\in\{0,1\}.
\]
Abbreviate
\[
 c=c_n,\quad t=t_n,\quad C=c_{n+1},\quad T=t_{n+1},
 \qquad d=c+n-1,
\]
and define
\begin{align}
 A(x)&=(T+x)x(c+x-1)
       -(t+x-1)(x+1)(C+x),\label{eq:circle-G-A}\\
 B(x)&=(n-x)(\overline c+n-x-1)(C+x)\nonumber\\
     &\hspace{9mm}-(n+1-x)(T+x)(\overline t+n-x-1).
       \label{eq:circle-G-B}
\end{align}
Then $A,B$ are nonzero polynomials of degree at most two, and
\begin{equation}\label{eq:circle-G-g}
 g(x)=\frac{(c+x-1)B(x)}
            {(\overline t+n-x-1)A(x)}
\end{equation}
satisfies \eqref{eq:circle-G-rational-shift}.  In particular,
\begin{equation}\label{eq:circle-G-three-shifts}
 \begin{gathered}
 m:=t-c\in\mathbb Z,\quad |m|\leq3,
 \quad\hbox{or}\quad
 j:=t+\overline c+n-2\in\mathbb Z,\quad |j|\leq3.
 \end{gathered}
\end{equation}
\end{lemma}

\begin{proof}
Write $\Phi_n=\sum_{k=0}^n a_kz^k$, $a_n=1$.
By Lemma~\ref{lem:circle-G-coefficient-structure} and
nonresonance, the nonzero coefficients are precisely
$a_r,\ldots,a_n$. If $r=1$, the constant-coefficient equation
in \eqref{eq:circle-G-coefficient-recurrence} gives $c=0$.

Since $\delta_{n+1}\neq0$, all coefficients of $\Phi_{n+1}$
are nonzero.  Its constant-coefficient equation also implies
$C\neq0$.

Set
\[
 U_k=a_{k-1},\qquad
 V_k=\delta_{n+1}\overline{a_{n-k}}.
\]
The coefficient of $z^k$ in the Szeg\H{o} recurrence is $U_k+V_k$.
For $r+1\leq k\leq n-2$, all of $U_k,V_k,U_{k+1},V_{k+1}$
are nonzero, and their successive ratios are
\begin{align*}
 R(k)&=-\frac{(n+1-k)(t+k-1)}{k(c+k-1)},\\
 S(k)&=-\frac{(n-k)(\overline c+n-k-1)}
                 {(k+1)(\overline t+n-k-1)},\\
 L(k)&=-\frac{(n+1-k)(T+k)}{(k+1)(C+k)}.
\end{align*}
Thus
\begin{equation}\label{eq:circle-G-coefficient-sum}
 [R(k)-L(k)]U_k=[L(k)-S(k)]V_k.
\end{equation}
Using this relation at $k$ and $k+1$, and eliminating $U_k,V_k$
without dividing by $A$ or $B$, yields $H(k)=0$, where
\begin{equation}\label{eq:circle-G-H}
 \begin{split}
 H(x)={}&(c+x)(\overline c+n-x-1)B(x+1)A(x)\\
 &-(t+x-1)(\overline t+n-x-2)B(x)A(x+1).
 \end{split}
\end{equation}
This holds for $r+1\leq k\leq n-3$; both indices used in
\eqref{eq:circle-G-coefficient-sum} are then within the range above.  Every denominator cleared
in this calculation is nonzero in that range, by nonresonance
and the nonvanishing coefficient intervals.  The cubic terms in
$A$ and $B$ cancel, and the degree-six terms in $H$ cancel.
Hence $\deg H\leq5$. The indicated range contains
$n-r-3\geq6$ distinct integer indices, so $H\equiv0$.

If $A\equiv0$, equation \eqref{eq:circle-G-coefficient-sum} forces
$B\equiv0$, and conversely.  Comparing coefficients in these two
identities, using $C\neq0$, gives
\[
 t=1,\qquad
 \left\{
 \begin{array}{ll}
 c=2,&C=T,\\
 c=-n,&T=1,\ C=-n-1.
 \end{array}\right.
\]
In these two cases the monic constant term of $\Phi_n$ has
modulus $n+1$ and $1$, respectively.  Both contradict positivity.
Thus $A$ and $B$ are nonzero.

The identity $H\equiv0$ is equivalent to
\eqref{eq:circle-G-rational-shift} for the function
\eqref{eq:circle-G-g}.  Its numerator and denominator have degree
at most three.  Lemma~\ref{lem:circle-G-rational-shift} now gives
\eqref{eq:circle-G-three-shifts}.
\end{proof}

\begin{lemma}[Local Gauss rigidity]
\label{lem:circle-G-rigidity}
Under the hypotheses of Lemma~\ref{lem:circle-G-compatibility},
one has
\begin{equation}\label{eq:circle-G-rigidity}
 \delta_n\neq0,\qquad t_n=1-\overline{d_n}.
\end{equation}
\end{lemma}

\begin{proof}
Use the notation of Lemma~\ref{lem:circle-G-compatibility}.
Suppose first that $r=1$, so $c=0$.  If $m=t-c$ is one of the
integers in \eqref{eq:circle-G-three-shifts}, nonresonance permits
only $t=1,2,3$.  Writing $\Phi_n=zQ$, the polynomial $Q$ has
degree $n-1$ and satisfies a Gauss equation with parameters
$c_Q=2$ and $t_Q=t+1$.  Vandermonde's identity gives $Q(1)=0$,
because $(2-t_Q)_{n-1}=(1-t)_{n-1}=0$.  This contradicts the
disk location of OPUC zeros.

If instead $j=t+\overline c+n-2$ is one of the indicated integers,
then $t=j-n+2$.  Nonresonance leaves only $j=-3,-2$, giving
$t=-n-1,-n$.  But the coefficient of $z^{n-1}$ is then
\[
 s_n=-\frac{n(n-1)}{n-1+t},
\]
whose modulus is at least $n$, impossible for zeros strictly inside
the unit disk.  Hence $r=0$, and $\delta_n\neq0$.

Coefficient comparison now gives
\begin{equation}\label{eq:circle-G-delta-s}
 \delta_n=\frac{(-d)_n}{(t)_n},\qquad
 s_n=-\frac{nd}{n-1+t}.
\end{equation}
The disk location of the zeros implies
\begin{equation}\label{eq:circle-G-stability}
 |\delta_n|<1,\qquad |d|<|n-1+t|,\qquad |t|>|d|.
\end{equation}
For the last inequality, evaluation at $z=1$ first shows that
$t\ne d$, since $t\ne0$ and $\Phi_n(1)\ne0$. Then use
\[
 \frac{\Phi_n'(1)}{\Phi_n(1)}=\frac{nt}{t-d},\qquad
 \operatorname{Re}\frac{\Phi_n'(1)}{\Phi_n(1)}>\frac n2.
\]
The same reasoning at degree $n+1$ gives
\begin{equation}\label{eq:circle-G-next-stability}
 |T|>|C+n|.
\end{equation}

We eliminate all alternatives in
\eqref{eq:circle-G-three-shifts} except $j=0$.

\emph{Case $m=0,1,2,3$.}
Vandermonde's identity gives
\[
 \Phi_n(1)=(-1)^n\frac{(-m)_n}{(t)_n}=0,
\]
again contradicting the location of the zeros.

\emph{Case $m=-l$, $1\leq l\leq3$.}
Here $t=c-l$, and \eqref{eq:circle-G-stability} gives
$\operatorname{Re}d<l/2$.  A rational solution of
\eqref{eq:circle-G-rational-shift} is
\[
 g(x)=\kappa\frac{(x-\overline d)_l}{(x+c-l-1)_l},
 \qquad \kappa\neq0.
\]
This is the general nonzero rational solution up to the constant:
the quotient of two rational solutions is one-periodic and hence
constant.  Comparing with \eqref{eq:circle-G-g} gives
\begin{equation}\label{eq:circle-G-parallel-products}
 B(x)(x+c-l-1)_{l+1}
 =-\kappa A(x)(x-\overline d)_{l+1}.
\end{equation}
The two displayed Pochhammer polynomials have disjoint root sets.
Indeed, the real parts of the roots on the left are at least
$n-\operatorname{Re}d$, whereas those on the right are at most
$\operatorname{Re}d$, and $n>l$.
Consequently $A(1-c)=0$ and $B(\overline d)=0$.
Using \eqref{eq:circle-G-A} and \eqref{eq:circle-G-B}, these
equalities give
\begin{equation}\label{eq:circle-G-forbidden-next}
 C=c-1,\qquad T=-\overline d.
\end{equation}
For clarity, the evaluations used here are
\begin{align*}
 A(1-c)&=-(t-c)(2-c)(C+1-c),\\
 B(\overline d)&=-(n+1-\overline d)(T+\overline d)
                    (\overline t+n-\overline d-1).
\end{align*}
The factors canceled are nonzero: $t-c=-l$,
$\overline t+n-\overline d-1=-l$, and
$\operatorname{Re}d<l/2$ with $n\geq10$.
But $C+n=d$, so \eqref{eq:circle-G-forbidden-next}
contradicts \eqref{eq:circle-G-next-stability}.

\emph{Case $j=-l$, $1\leq l\leq3$.}
Now $t=1-\overline d-l$.  If $l=1$, equation
\eqref{eq:circle-G-delta-s} gives $|\delta_n|=1$.
If $l=2$ or $3$, the same equation and the middle inequality in
\eqref{eq:circle-G-stability} give
\begin{equation}\label{eq:circle-G-negative-cross}
 |\delta_n|=
 \prod_{h=1}^{l-1}\frac{|n-d-h|}{|d+h|},
 \qquad \operatorname{Re}d<\frac{n-l}{2}.
\end{equation}
For $l=2$, comparison of squared moduli immediately makes the
product greater than one.  For $l=3$, put
$x=\operatorname{Re}d+3/2<n/2$, $y=\operatorname{Im}d$, and
\[
 F(x)=\bigl((x-\tfrac12)^2+y^2\bigr)
       \bigl((x+\tfrac12)^2+y^2\bigr).
\]
Then
\[
 F(n-x)-F(x)=
 n(n-2x)\bigl((n-x)^2+x^2+2y^2-\tfrac12\bigr)>0.
\]
Thus the product in \eqref{eq:circle-G-negative-cross} is again
greater than one, contradicting $|\delta_n|<1$.

\emph{Case $j=1,2,3$.}
Here $t=1-\overline d+j$.  The last inequality in
\eqref{eq:circle-G-stability} gives
$\operatorname{Re}d<(j+1)/2$.  Solving the rational difference
equation gives
\[
 g(x)=\kappa\frac{(x-\overline d)_j}{(x+c-1-j)_j},
 \qquad\kappa\neq0,
\]
and hence
\begin{equation}\label{eq:circle-G-positive-cross}
 B(x)(x+c-j)_j=-\kappa A(x)(x-\overline d)_j.
\end{equation}
The two Pochhammer root sets are disjoint, since their real parts
are separated by $n-2\operatorname{Re}d>n-j-1>0$.
Thus $A(1-c)=B(\overline d)=0$.
Here $t-c=\overline t+n-\overline d-1
=n+j-2\operatorname{Re}d>n-1$, so the same evaluations again
give \eqref{eq:circle-G-forbidden-next}, contrary to
\eqref{eq:circle-G-next-stability}.

Only $j=0$ remains.  Since
$j=t+\overline c+n-2=t+\overline d-1$, this gives
$t=1-\overline d$, proving the lemma.
\end{proof}

\paragraph{Identification of the measure.}

The coefficient relation obtained in
Lemma~\ref{lem:circle-G-rigidity} leads to the circular
hypergeometric family studied
in~\cite{SriRanga2010,BorregoSriRanga2017}.
We first establish the properties of this family and
identify the corresponding measure under a fixed-parameter
coefficient law imposed only for $n\ge N_0$. The subsequent
rigidity theorem derives this law from the general hypotheses,
allowing the stated freedom at degree one.

\begin{theorem}[The positive circular hypergeometric family]
\label{thm:circle-Jacobi-branch}
For $q\in\mathbb C$ with $\operatorname{Re}q>-\tfrac12$,
define the positive probability measure
\begin{equation}
\label{eq:circle-G-circular-measure}
    d\nu_q(e^{i\varphi})=
    C_q e^{(\pi-\varphi)\operatorname{Im}q}
    \left(2\sin\frac{\varphi}{2}\right)^{2\operatorname{Re}q}
    \frac{d\varphi}{2\pi},
    \qquad 0<\varphi<2\pi,
\end{equation}
where $C_q>0$ is the normalizing constant.
Its monic OPUC are
\begin{equation}
\label{eq:circle-G-regular-polynomial}
    P_n^{(q)}(z)=
    \frac{(q+\overline q+1)_n}{(q+1)_n}
    {}_2F_1(-n,q+1;q+\overline q+1;1-z).
\end{equation}
For every $n\ge1$, they satisfy \eqref{eq:circle-G} with
\begin{equation}
\label{eq:circle-G-canonical-coefficients}
    b_n=q+2-n,\qquad c_n=1-n-\overline q,
    \qquad \lambda_n=n(q+1).
\end{equation}
Their constant terms and Verblunsky coefficients are
\begin{equation}
\label{eq:circle-G-circular-Verblunsky}
    P_n^{(q)}(0)=\frac{(\overline q)_n}{(q+1)_n},
    \qquad
    \alpha_n^{(q)}
    =-\frac{(q)_{n+1}}{(\overline q+1)_{n+1}}.
\end{equation}

Furthermore, the Laurent moment functional
$u_q(f)=\int_{\mathbb T}f\,d\nu_q$ satisfies the
Pearson identity
\begin{equation}
\label{eq:circle-Jacobi-Pearson}
    D(\sigma u_q)=\tau_q u_q,\qquad
    B(z,n)=\tau_q(z)+(1-n)\frac{\sigma(z)}{z},
\end{equation}
where
\[
    \sigma(z)=z(1-z),\qquad
    \tau_q(z)=-\overline q-(q+1)z,
\]
and $Du_q(f)=-u_q(f')$ for Laurent polynomials.
The corresponding angular Pearson identity is
\begin{equation}
\label{eq:circle-G-candidate-Pearson}
    D_\varphi[(1-e^{i\varphi})\nu_q]
    =-i[\overline q+(q+1)e^{i\varphi}]\nu_q
\end{equation}
in distributions on the circle.

Conversely, let $\mu$ be a positive probability measure
of infinite support on $\mathbb T$ whose monic OPUC
satisfy \eqref{eq:circle-G} for every $n\ge N_0$ with
\eqref{eq:circle-G-canonical-coefficients} for some fixed
$q\in\mathbb C$. Then necessarily
$\operatorname{Re}q>-\tfrac12$, $\mu=\nu_q$, and
$\Phi_n=P_n^{(q)}$ for every $n\ge0$.
In particular, there is no additional atom at $1$.
The choice $q=0$ gives the normalized Lebesgue measure.
\end{theorem}

\begin{proof}
First suppose that $\operatorname{Re}q>-\tfrac12$.
Write $q=a+i\eta$ and put
$\gamma=q+\overline q+1=2a+1>0$.
The weight in \eqref{eq:circle-G-circular-measure}
is positive and integrable, since $2a>-1$.

Define
\[
    W(z)=z^{-\overline q}(1-z)^{2a},
\]
with branches chosen along the unit circle cut at $1$.
Indeed, taking $\arg z=\varphi$ and
$\arg(1-z)=(\varphi-\pi)/2$, we obtain
\[
    W(e^{i\varphi})
    =e^{-\pi(\eta+ia)}
      e^{(\pi-\varphi)\eta}
      \left(2\sin\frac{\varphi}{2}\right)^{2a}.
\]
Thus, up to a nonzero constant, the boundary values of
$W$ coincide with the density in
\eqref{eq:circle-G-circular-measure}.

Consider the  Rodrigues expression
\begin{equation}
\label{eq:circle-Jacobi-Rodrigues}
    P_n(z)=\frac{(-1)^n}{(q+1)_n}
         \frac{z^n}{W(z)}
         D^n\bigl[(1-z)^nW(z)\bigr],
\end{equation}
where $D=d/dz$.
By Leibniz's rule,
\[
    P_n(z)=\frac{1}{(q+1)_n}
    \sum_{j=0}^n\binom nj
    (\overline q)_j(\gamma+j)_{n-j}
    z^{n-j}(1-z)^j.
\]
Since $(\gamma+j)_{n-j}=(\gamma)_n/(\gamma)_j$,
the terminating Pfaff transformation gives
\[
\begin{aligned}
    P_n(z)
    &=
    \frac{(\gamma)_n}{(q+1)_n}
    z^n{}_2F_1
       \left(-n,\overline q;\gamma;\frac{z-1}{z}\right)\\
    &=
    \frac{(\gamma)_n}{(q+1)_n}
    {}_2F_1(-n,q+1;\gamma;1-z).
\end{aligned}
\]
Consequently, $P_n=P_n^{(q)}$ as given by
\eqref{eq:circle-G-regular-polynomial}.
The coefficient of $z^n$ in this expression is one.
The Gauss hypergeometric equation, under the substitution
$x=1-z$, gives \eqref{eq:circle-G} with
\eqref{eq:circle-G-canonical-coefficients}.

For $0\leq k<n$, the contour orthogonality integral
of $P_n$ against $z^{-k}$ is, apart from a nonzero constant,
\[
    \int_{\mathbb T\setminus\{1\}}
       z^{n-k-1}D^n[(1-z)^nW(z)]\,dz.
\]
Integration by parts $n$ times gives zero because
$D^nz^{n-k-1}=0$. All boundary terms vanish, since
\[
    D^j[(1-z)^nW(z)]
    =O(|1-z|^{n+2a-j}),
    \qquad 0\leq j\leq n-1,
\]
and $n+2a-j\geq2a+1>0$ at both endpoints of the
cut contour.
Hence the polynomials $P_n^{(q)}$ are precisely
the monic OPUC of $\nu_q$.

The constant term follows directly from the term $j=n$
in the Leibniz expansion, and
$\alpha_n^{(q)}=-\overline{P_{n+1}^{(q)}(0)}$
gives \eqref{eq:circle-G-circular-Verblunsky}.

To establish the Pearson identities, observe that
\[
    zD[(1-z)W(z)]=\tau_q(z)W(z).
\]
Since $d/d\varphi=izD$ along the circle, this gives
\eqref{eq:circle-G-candidate-Pearson} away from $1$.
Moreover, $(1-e^{i\varphi})W(e^{i\varphi})$ tends to zero
at both endpoints of the cut contour, because $2a+1>0$.
Integration by parts therefore establishes the same
identity globally in distributions.

Testing the angular identity against an arbitrary Laurent
polynomial $f$, and using
\[
    \frac{d}{d\varphi}f(e^{i\varphi})
    =izf'(z),
\]
yields
\[
    -u_q(\sigma f')=u_q(\tau_q f),
\]
which is precisely the Laurent Pearson identity
\eqref{eq:circle-Jacobi-Pearson}.
The expression for $B(z,n)$ follows directly from
\eqref{eq:circle-G-canonical-coefficients}, exhibiting
its degree dependence in terms of the fixed Pearson
pair $(\sigma,\tau_q)$.

In particular, writing
$M_k=\int_{\mathbb T}z^k\,d\nu_q(z)$, the Pearson
identity gives
\[
    (k-\overline q)M_k=(k+q+1)M_{k+1},
\]
and hence
\[
    M_k=\frac{(-\overline q)_k}{(q+1)_k},
    \qquad k\geq0.
\]
The negative moments are their conjugates.

For the converse, suppose that the monic OPUC of $\mu$
satisfy \eqref{eq:circle-G} for every $n\ge N_0$ with
\eqref{eq:circle-G-canonical-coefficients} for a fixed
$q\in\mathbb C$.
Put $a=\operatorname{Re}q$. Choose any $n\ge N_0$.
Evaluation of the differential equation at $z=1$ gives
\[
    (2a+1)\Phi_n'(1)=n(q+1)\Phi_n(1).
\]
Since $\Phi_n(1)\neq0$ and $q+1\neq0$ when $a=-1/2$,
we have $2a+1\neq0$. Because all zeros $\xi$ of $\Phi_n$
satisfy $|\xi|<1$, each contributes
$\operatorname{Re}(1/(1-\xi))>1/2$ to the real part of its
logarithmic derivative at $1$. Thus
\[
    \frac n2
    <\operatorname{Re}\frac{\Phi_n'(1)}{\Phi_n(1)}
    =\frac{n(a+1)}{2a+1}
    =\frac n2+\frac{n}{2(2a+1)}.
\]
Consequently, $a>-1/2$.

Under \eqref{eq:circle-G-canonical-coefficients},
the recurrence \eqref{eq:circle-G-coefficient-recurrence}
determines the coefficients of a monic polynomial solution
uniquely, since $t_n=q+1$ and
\[
    q+1+k\neq0,\qquad 0\leq k<n.
\]
Consequently, $\Phi_n=P_n^{(q)}$ for every $n\ge N_0$.
The inverse Szeg\H{o} recurrence applied to both sequences
gives the same equality in every earlier degree.
Uniqueness of the probability measure determined by its
full monic OPUC sequence gives $\mu=\nu_q$.

This also excludes any additional atom at $1$.
Independently, since $\sigma(1)=0$ whereas
$\tau_q(1)=-(2a+1)\neq0$, adding a nonzero mass at $1$
would violate the Pearson identity.
Finally, $q=0$ gives the monomials and the normalized Lebesgue measure.
\end{proof}

\begin{theorem}[Gauss-to-Pearson rigidity in the infinite branch]
\label{thm:circle-infinite-G}
Let $\mu$ be a positive probability measure of infinite support
on $\mathbb{T}$, and let $\Phi_n$ be its monic OPUC. Assume
\eqref{eq:circle-G} for every $n\ge N_0$. If infinitely many
Verblunsky coefficients are nonzero,
then there exists a unique
\[
 q\in\mathbb C\setminus\{0\},\qquad
 \operatorname{Re}q>-\tfrac12,
\]
such that $\mu=\nu_q$, with $\nu_q$ given by
\eqref{eq:circle-G-circular-measure}.  In particular,
\begin{equation}\label{eq:circle-G-Pearson}
 \displaystyle
 D_\varphi[(1-e^{i\varphi})\mu]
 =-i[\overline q+(q+1)e^{i\varphi}]\mu.
\end{equation}
Thus the common Pearson constants are
\[
 d=-\overline q,\qquad t=q+1.
\]
In all degrees the polynomials and Verblunsky coefficients are
given by
\eqref{eq:circle-G-regular-polynomial} and
\eqref{eq:circle-G-circular-Verblunsky}.  Every Verblunsky
coefficient is nonzero. The only admissible differential
coefficients in degrees $n\ge2$ are
\[
 b_n=q+2-n,\qquad c_n=1-n-\overline q\qquad(n\geq2).
\]
In particular, these coincide with the prescribed coefficients
for $n\ge\max\{2,N_0\}$ and provide an extension to all
earlier degrees $n\ge2$. At degree one, $b_1\in\mathbb C$
may be chosen arbitrarily and
$c_1=-b_1\overline q/(q+1)$; the canonical choice in
\eqref{eq:circle-G-canonical-coefficients} is admissible there too.
Conversely, every $\nu_q$ with $q\neq0$ and
$\operatorname{Re}q>-1/2$ satisfies \eqref{eq:circle-G}
in every degree $n\ge1$ and has infinitely many nonzero
Verblunsky coefficients.
\end{theorem}

\begin{proof}
Since infinitely many Verblunsky coefficients are nonzero,
the measure is not the normalized Lebesgue measure.
Lemma~\ref{lem:circle-G-automatic-nonresonance}
guarantees nonresonance at every degree $n\ge N_1$.
Lemma~\ref{lem:circle-G-zero-tail} excludes consecutive zero
constant terms whenever the first degree is at least
$\max\{2,N_0\}$.
If $\delta_n=0$ for some $n\ge N_2$, then
$n-1\ge\max\{2,N_0\}$, and hence
$\delta_{n-1}\neq0$ and $\delta_{n+1}\neq0$.
The Szeg\H{o} recurrence gives $\Phi_n=z\Phi_{n-1}$, so
$\operatorname{ord}_0\Phi_n=1$.
As $n\ge10$ and both degrees $n,n+1$ are in the nonresonant
range, Lemma~\ref{lem:circle-G-rigidity} applies and implies
$\delta_n\neq0$, a contradiction. Therefore
\[
 \delta_n\neq0\qquad(n\ge N_2).
\]
Applying the same rigidity lemma at every such degree gives
\[
 t_n=1-\overline{d_n},\qquad n\ge N_2.
\]
Put $q_n=-\overline{d_n}$.  The strict inequality
$|t_n|>|d_n|$ from \eqref{eq:circle-G-stability} gives
$\operatorname{Re}q_n>-1/2$.
At this degree the differential equation has exactly the
canonical coefficients with parameter $q_n$.  Nonresonance
ensures uniqueness of the monic polynomial solution, so
$\Phi_n=P_n^{(q_n)}$.

Inverse Szeg\H{o} recurrence uniquely determines every earlier
polynomial from $\Phi_n$.  Applying it also to the positive
family in Theorem~\ref{thm:circle-Jacobi-branch} gives
\[
 \delta_1=\frac{\overline{q_n}}{q_n+1}.
\]
This determines $q_n$ uniquely:
\begin{equation}\label{eq:circle-G-q-from-first-root}
 q_n=\frac{|\delta_1|^2+\overline{\delta_1}}
           {1-|\delta_1|^2}.
\end{equation}
Thus all $q_n$ equal one constant $q$. For any $m\ge0$, take
$n\ge\max\{N_2,m\}$. The inverse recurrence then gives
$\Phi_m=P_m^{(q)}$.  The complete monic OPUC
sequence determines its probability measure, so $\mu=\nu_q$.
Equation \eqref{eq:circle-G-Pearson} follows from
\eqref{eq:circle-G-candidate-Pearson}.  The parameter is nonzero,
because $q=0$ would give monomials.

For $q\neq0$ the constant terms in
\eqref{eq:circle-G-circular-Verblunsky} never vanish, and
\[
 |\alpha_n|^2=\frac{|q|^2}{|q+n+1|^2}<1.
\]
By Theorem~\ref{thm:circle-Jacobi-branch}, the recovered
polynomials satisfy the canonical Gauss equation in all degrees.
For $n\geq2$, the admissible coefficients $b_n,c_n$ are unique. Indeed,
subtracting two possible Gauss equations gives
\[
 (\Delta c-\Delta b\,z)\Phi_n'+n\Delta b\,\Phi_n=0.
\]
The polynomial has at least two distinct nonzero zeros, all
ordinary points of the equation.  Evaluation at these zeros
forces $\Delta b=\Delta c=0$.
This proves the asserted coefficient law in every prescribed
degree $n\ge2$. At degree one direct substitution gives the
stated freedom.

Conversely, Theorem~\ref{thm:circle-Jacobi-branch} gives the
Gauss equation in all positive degrees, positivity, and
infinite support. Since
$\operatorname{Re}(q+1)>1/2$, the canonical coefficients satisfy
nonresonance, and $q\neq0$ makes all Verblunsky coefficients
nonzero.
\end{proof}

\paragraph{Conclusion of the Gauss classification.}

Combining the finite and infinite Verblunsky cases gives the
following complete classification.

\begin{corollary}[Classification of the Gauss case]
\label{cor:circle-G-complete}
The positive probability measures of infinite support on
$\mathbb T$ whose monic OPUC satisfy \eqref{eq:circle-G}
for every $n\ge N_0$, for some integer $N_0\ge1$, are exactly
the following:
\begin{enumerate}
\item the Poisson measures of Theorem~\ref{thm:circle-finite-G},
including the normalized Lebesgue measure;
\item the circular hypergeometric measures
\eqref{eq:circle-G-circular-measure} with $q\neq0$ and
$\operatorname{Re}q>-1/2$.
\end{enumerate}
The first class has finitely many nonzero Verblunsky coefficients;
in the second class every Verblunsky coefficient is nonzero.
For both non-Lebesgue families, the admissible differential
coefficients are necessarily nonresonant in every degree
$n\ge2$.
At degree one, the resonant choice $b_1=c_1=0$ is admissible.
For Lebesgue, $c_n=1-n$ and $b_n$ is arbitrary in all degrees.
In every case the initially prescribed differential equations
can be extended to all degrees $n\ge1$ with the same leading
coefficient. Thus imposing \eqref{eq:circle-G} only for
sufficiently large degrees yields exactly the same
classification as imposing it in every positive degree.
In particular, no measure in this classification has a nonzero
atomic or singular continuous part.
\end{corollary}

\begin{proof}
Apply Theorem~\ref{thm:circle-finite-G} when only finitely many
Verblunsky coefficients are nonzero, and
Theorem~\ref{thm:circle-infinite-G} otherwise. The converse
statements in those theorems show that every listed measure
admits \eqref{eq:circle-G} in all positive degrees.
The restrictions on the coefficients proved there give the
asserted extension.
\end{proof}

\subsection{Classification of the double-zero case $\theta(z)=z^2$}
\label{sec:circle-double-zero}

We consider equation $\mathrm{(E)}$, including the solvable
resonant cases described in
Proposition~\ref{prop:circle-polynomial-solvability}.
The equation is of Euler type when $c_n=0$ and of generalized
Bessel type otherwise. As will be shown, resonance is essential
for the occurrence of non-Lebesgue measures in this branch.
Throughout this subsection, the orthogonality circle is $\mathbb T$
and the double leading zero is its center. We allow $\mathrm{(E)}$
to be imposed only in sufficiently large degrees.

Put
\[
    t_n=n-1-b_n,\qquad
    \delta_n=\Phi_n(0)=-\overline{\alpha_{n-1}}.
\]
Whenever $\Phi_n$ satisfies $\mathrm{(E)}$, leading-coefficient
comparison gives $\lambda_n=-nt_n$. Writing
$\Phi_n(z)=\sum_{j=0}^n a_{n,j}z^j$, $a_{n,n}=1$, we specialize
\eqref{eq:circle-polynomial-solvability-recurrence} to the present
case:
\begin{equation}
\label{eq:circle-E-coefficient-recurrence}
    c_n(j+1)a_{n,j+1}=(n-j)(j+t_n)a_{n,j},
    \qquad 0\le j<n.
\end{equation}
If $c_n\ne0$, downward recursion from $a_{n,n}=1$ shows that
$j+t_n\ne0$ and $a_{n,j}\ne0$ for every $0\le j<n$.
Thus the equation is nonresonant and every coefficient of
$\Phi_n$ is nonzero.
If $c_n=0$, the same recurrence shows that $\Phi_n=z^n$
or $\Phi_n=z^n+d z^r$, with $d\ne0$, $0\le r<n$, and
$t_n=-r$ in the latter case. Moreover, $|d|<1$, since the
nonzero zeros of this binomial have modulus $|d|^{1/(n-r)}$.

\subsubsection{Compatibility with the Szeg\H{o} recurrence}

We first show that the occurrence of $c_N=0$ at a degree
$N\ge3$ restricts the entire OPUC sequence, provided that
$\mathrm{(E)}$ holds in every subsequent degree. We then prove
that such a degree necessarily exists whenever $\mathrm{(E)}$
holds in all sufficiently large degrees. This will yield a strengthening of
Theorem~\ref{thm:circle-main-classification} in the double-zero
case: the same classification holds when $\mathrm{(E)}$ is
imposed only in sufficiently large degrees.

\begin{lemma}[Propagation from a degree with $c_N=0$]
\label{lem:circle-E-binomial-propagation}
Suppose that, for some $N\ge3$, $\mathrm{(E)}$ holds for
every $n\ge N$ and that $c_N=0$. Then at most one Verblunsky
coefficient is nonzero.
\end{lemma}

\begin{proof}
Suppose first that at a degree $n\ge N$ one has
\[
    \Phi_n=z^n+d z^r,
    \qquad d\ne0,\quad 0\le r<n.
\]
Put $\varepsilon=\delta_{n+1}$. If $\varepsilon\ne0$,
the Szeg\H{o} recurrence gives
\begin{equation}
\label{eq:circle-E-next-binomial}
    \Phi_{n+1}=z^{n+1}+d z^{r+1}
                    +\varepsilon
                    +\varepsilon\overline d\,z^{n-r}.
\end{equation}
This polynomial has at most four nonzero coefficients.
Since $c_{n+1}\ne0$ would make all $n+2\ge5$ coefficients
nonzero, we must have $c_{n+1}=0$. Its nonzero constant
term then forces $\Phi_{n+1}=z^{n+1}+\varepsilon$.
The two intermediate terms in
\eqref{eq:circle-E-next-binomial} cannot cancel unless
$r+1=n-r$ and $d+\varepsilon\overline d=0$.
The latter equality gives $|\varepsilon|=1$, contrary to
positivity. Therefore $\delta_{n+1}=0$.

It follows that $\Phi_{n+1}=z\Phi_n$ is again a nonmonomial
binomial. Since the differential equation holds in every
subsequent degree, induction gives
$\Phi_{n+k}=z^k\Phi_n$ for every $k\ge0$.
Set $m=n-r$. Applying the inverse Szeg\H{o} recurrence $r$
times first gives $\Phi_m=z^m+d$. Since $|d|<1$, it then gives
\[
    \Phi_{m-1}
       =\frac{\Phi_m-d\Phi_m^*}{z(1-|d|^2)}=z^{m-1}.
\]
All earlier polynomials are monomials. Thus only
$\delta_m=d$ is nonzero. Notice that this backward argument
uses no differential equation in degrees below $N$.

Now use $c_N=0$. If $\Phi_N$ is nonmonomial, the preceding
argument applies. If $\Phi_N=z^N$, all earlier polynomials
are monomials. Either every later constant term vanishes,
or there is a first $M>N$ with $\delta_M\ne0$.
In the latter case $\Phi_M=z^M+\delta_M$, with $M\ge4$,
so the preceding binomial argument applies at degree $M$.
In every case at most one Verblunsky coefficient is nonzero.
\end{proof}

\begin{lemma}[Exclusion of consecutive nonzero $c_n$]
\label{lem:circle-E-consecutive-dense}
Let $n\ge10$. Suppose that $\Phi_n,\Phi_{n+1}$ satisfy
$\mathrm{(E)}$ at their respective degrees and the Szeg\H{o}
recurrence. Then
\[
    c_nc_{n+1}=0.
\]
\end{lemma}

\begin{proof}
Suppose instead that both parameters are nonzero. Abbreviate
\[
    c=c_n,\quad C=c_{n+1},\quad
    t=t_n,\quad T=t_{n+1},\qquad a_j=a_{n,j}.
\]
Both equations are nonresonant, and
\eqref{eq:circle-E-coefficient-recurrence} implies that all
coefficients of both polynomials are nonzero. In particular,
$\delta_{n+1}\ne0$.
For the coefficient of $z^k$ in
$\Phi_{n+1}=z\Phi_n+\delta_{n+1}\Phi_n^*$, put
\[
    U_k=a_{k-1},\qquad
    V_k=\delta_{n+1}\overline{a_{n-k}}.
\]
For $1\le k\le n-1$, define the successive coefficient ratios
\begin{align*}
    R(k)&=\frac{(n+1-k)(k+t-1)}{ck},\\
    S(k)&=\frac{\overline c(n-k)}
                   {(k+1)(\overline t+n-k-1)},\\
    L(k)&=\frac{(n+1-k)(k+T)}{C(k+1)}.
\end{align*}
Here $U_{k+1}=R(k)U_k$, $V_{k+1}=S(k)V_k$, and
$L(k)$ is the successive ratio for the coefficients of
$\Phi_{n+1}$. Therefore
\begin{equation}
\label{eq:circle-E-coefficient-compatibility}
    [R(k)-L(k)]U_k=[L(k)-S(k)]V_k.
\end{equation}
Define
\begin{align}
\label{eq:circle-E-A}
    A(x)&=C(x+1)(x+t-1)-cx(x+T),\\
\label{eq:circle-E-D}
    D(x)&=(n+1-x)(x+T)(\overline t+n-x-1)
                  -C\overline c(n-x).
\end{align}
Then
\[
    R(k)-L(k)=\frac{(n+1-k)A(k)}{cCk(k+1)},
    \qquad
    L(k)-S(k)=\frac{D(k)}
                     {C(k+1)(\overline t+n-k-1)}.
\]
The degree of $A$ is at most two, and $D$ is monic of
degree three. If $A\equiv0$, then
\eqref{eq:circle-E-coefficient-compatibility}, with $V_k\ne0$,
forces $D(k)=0$ for $1\le k\le n-1$. This contradicts the
degree of $D$. Thus $A\not\equiv0$; put $d=\deg A\le2$.

Use \eqref{eq:circle-E-coefficient-compatibility} at $k$
and $k+1$. Eliminating $U_k,V_k$ without division by $A$ or $D$
gives
\[
\begin{split}
    &[L(k)-S(k)][R(k+1)-L(k+1)]R(k)\\
    &\qquad-[R(k)-L(k)][L(k+1)-S(k+1)]S(k)=0.
\end{split}
\]
After clearing nonzero denominators and cancelling the nonzero
factors $n-k$ and $n+1-k$, this yields
\begin{equation}
\label{eq:circle-E-H-zeros}
    H(k)=0,\qquad 1\le k\le n-2,
\end{equation}
where
\begin{equation}
\label{eq:circle-E-H}
\begin{split}
    H(x)={}&|c|^2D(x+1)A(x)\\
            &-(x+t-1)(\overline t+n-x-2)D(x)A(x+1).
\end{split}
\end{equation}
All denominators involving $t$ are nonzero in this index range
by nonresonance; the remaining denominators are nonzero
because $cC\ne0$ and $k\ge1$.
The first term in $H$ has degree $d+3$, whereas the second
has degree $d+5$ with nonzero leading coefficient.
Hence $H\not\equiv0$ and $\deg H=d+5\le7$.
But \eqref{eq:circle-E-H-zeros} gives $n-2\ge8$ distinct
zeros, a contradiction.
\end{proof}

\subsubsection{Classification and admissible parameters}

\begin{theorem}[Complete positive OPUC classification for $\theta=z^2$]
\label{thm:circle-E-classification}
Let $\mu$ be a positive probability measure of infinite support
on $\mathbb T$, and let $\Phi_n$ be its monic orthogonal
polynomials. There exist $N_0\ge1$ and complex sequences
$b_n,c_n$ such that $\Phi_n$ satisfies $\mathrm{(E)}$ for every
$n\ge N_0$ if and only if one of the following holds.
\begin{enumerate}
\item The measure is the normalized Lebesgue measure and
$\Phi_n=z^n$. In every degree in which $\mathrm{(E)}$ is
imposed, its admissible differential parameters are $c_n=0$,
with $b_n$ arbitrary.

\item There are a unique integer $m\ge1$ and a unique
$a\in\mathbb D\setminus\{0\}$ such that
\begin{equation}
\label{eq:circle-E-measure}
    d\mu_{m,a}(e^{i\varphi})
       =\frac{1-|a|^2}{|e^{im\varphi}-a|^2}
          \frac{d\varphi}{2\pi}.
\end{equation}
The monic orthogonal polynomials and Verblunsky coefficients are
\begin{equation}
\label{eq:circle-E-polynomials}
\begin{gathered}
    \Phi_n(z)=
    \begin{cases}
        z^n,&0\le n<m,\\
        z^{n-m}(z^m-a),&n\ge m,
    \end{cases}\\
    \alpha_j=
    \begin{cases}
        \overline a,&j=m-1,\\
        0,&j\ne m-1.
    \end{cases}
\end{gathered}
\end{equation}
One admissible choice in every positive degree is
\begin{equation}
\label{eq:circle-E-canonical-parameters}
    c_n=0,\qquad b_n=2n-m-1,
    \qquad \lambda_n=n(n-m).
\end{equation}
Equivalently, the polynomials satisfy the Euler equations
\begin{equation}
\label{eq:circle-E-canonical-equation}
    z^2\Phi_n''-(2n-m-1)z\Phi_n'
        +n(n-m)\Phi_n=0.
\end{equation}
\end{enumerate}
In the non-Lebesgue case, all admissible differential parameters
in any positive degree in which $\mathrm{(E)}$ is imposed are
characterized by
\begin{equation}
\label{eq:circle-E-admissible-parameters}
\begin{cases}
    c_n=0,\quad b_n\in\mathbb C,
        &1\le n<m,\\
    c_1=ab_1,\quad b_1\in\mathbb C,
        &n=m=1,\\
    c_n=0,\quad b_n=2n-m-1,
        &n\ge\max\{m,2\}.
\end{cases}
\end{equation}
In both cases, the spectral coefficient in every imposed
degree is necessarily $\lambda_n=-n(n-1-b_n)$.
Every measure listed admits a realization of $\mathrm{(E)}$
in all positive degrees. In particular, all these measures
are absolutely continuous, with no atomic or singular
continuous part.
\end{theorem}

\begin{proof}
Suppose that $\mathrm{(E)}$ holds for every $n\ge N_0$.
Put $N_*=\max\{10,N_0\}$. By
Lemma~\ref{lem:circle-E-consecutive-dense},
$c_{N_*}c_{N_*+1}=0$. Hence there is a degree $N\ge N_*$
with $c_N=0$ for which the hypotheses of
Lemma~\ref{lem:circle-E-binomial-propagation} are satisfied.
It follows that at most one Verblunsky coefficient is nonzero.
If all vanish, the sequence consists of monomials and the
measure is the normalized Lebesgue measure.
Otherwise, let $\delta_m=-a\ne0$ be the unique nonzero
constant term. Positivity gives $0<|a|<1$, and the
Szeg\H{o} recurrence gives
\eqref{eq:circle-E-polynomials} in all degrees.

It remains to identify the measure and verify the converse.
The density in \eqref{eq:circle-E-measure} is positive
and normalized. Expanding its Poisson kernel gives
\begin{equation}
\label{eq:circle-E-moments}
    \int_{\mathbb T}z^k\,d\mu_{m,a}(z)=
    \begin{cases}
        a^{k/m},&m\mid k,\\
        0,&m\nmid k,
    \end{cases}
    \qquad k\ge0.
\end{equation}
Thus $z^n$ is orthogonal to all lower monomials when $n<m$.
For $n\ge m$ and $0\le j<n$, use
\[
    |z^m-a|^2
       =\frac{(z^m-a)(1-\overline a z^m)}{z^m}
       \qquad(z\in\mathbb T)
\]
to obtain
\[
    \int_{\mathbb T}\Phi_n(z)z^{-j}\,d\mu_{m,a}(z)
       =(1-|a|^2)\int_0^{2\pi}
           \frac{e^{i(n-j)\varphi}}
                {1-\overline a e^{im\varphi}}
           \frac{d\varphi}{2\pi}=0.
\]
The integrand has an absolutely convergent Fourier series
with only strictly positive powers. Hence the polynomials
in \eqref{eq:circle-E-polynomials} are the monic OPUC of
$\mu_{m,a}$. Uniqueness of the probability measure determined
by the complete monic sequence identifies the original $\mu$.
The single nonzero Verblunsky coefficient determines $m$ and
$a$ uniquely.

With $\mathcal D=z\,d/dz$, the differential operator in
\eqref{eq:circle-E-canonical-equation} factors as
\[
    (\mathcal D-n)(\mathcal D-n+m).
\]
It annihilates $z^n$ and, when $n\ge m$, $z^{n-m}$.
Thus \eqref{eq:circle-E-canonical-equation} holds in every
positive degree, proving the converse for $\mu_{m,a}$.
For Lebesgue, direct substitution gives $c_n=0$ and leaves
$b_n$ arbitrary.

To verify the full parameter description
\eqref{eq:circle-E-admissible-parameters}, consider any
positive degree in which $\mathrm{(E)}$ is imposed.
In the non-Lebesgue case, the polynomials in
\eqref{eq:circle-E-polynomials} have at most two nonzero
coefficients. Thus $c_n=0$ for $n\ge2$, since otherwise all
coefficients would be nonzero. If also $n\ge m$, the nonzero
coefficient at $z^{n-m}$ forces $t_n=m-n$, or
$b_n=2n-m-1$, by \eqref{eq:circle-E-coefficient-recurrence}.
For $n<m$, substitution of the monomial gives $c_n=0$
with $b_n$ arbitrary. Finally, when $n=m=1$, substitution
of $\Phi_1=z-a$ gives $c_1=ab_1$. Conversely, each of these
choices satisfies $\mathrm{(E)}$ in the indicated degree,
with $\lambda_n=-n(n-1-b_n)$. This proves the complete
parameter description.
\end{proof}

For the non-Lebesgue family $\mu_{m,a}$, the differential equation
is necessarily resonant in every imposed degree
$n\ge\max\{m,2\}$: its two Euler exponents are $n$ and $n-m$.
In degrees $1\le n<m$ in which the equation is imposed, it
may be resonant or nonresonant, with the monomial solution
selected in either case. If the equation is imposed at $n=m=1$,
it is resonant precisely when $b_1=0$.
For the normalized Lebesgue measure, the monomials solve $\mathrm{(E)}$ with
$c_n=0$ and arbitrary $b_n$ in every positive degree.
These possibilities cover all mixtures of resonant and
nonresonant degrees.

\begin{corollary}[The nonresonant double-zero case]
\label{cor:circle-E-nonresonant}
Under the hypotheses of
Theorem~\ref{thm:circle-E-classification}, if the monic
polynomial solution of $\mathrm{(E)}$ is unique in infinitely
many degrees in which the equation is imposed, then $\mu$
is the normalized Lebesgue measure. Conversely, Lebesgue measure admits a realization
that is nonresonant in every positive degree.
\end{corollary}

\begin{proof}
For every non-Lebesgue family,
\eqref{eq:circle-E-admissible-parameters} gives
$b_n=2n-m-1\in\{n-1,\ldots,2n-2\}$ for all sufficiently
large imposed degrees. By
Proposition~\ref{prop:circle-polynomial-solvability}, the monic
polynomial solution is nonunique at every such degree.
For the normalized Lebesgue measure, take $c_n=0$ and $b_n=-1$ for every $n\ge1$;
then $t_n=n$ is nonresonant.
\end{proof}

The preceding subsections exhaust the four canonical
configurations of the leading coefficient.
We are now in a position to combine their conclusions
and establish the equivalence asserted in
Theorem~\ref{thm:circle-main-classification}.

\begin{proof}[Proof of Theorem~\ref{thm:circle-main-classification}]
Comparison of the leading coefficients in
\eqref{eq:circle-main-DE} gives
\eqref{eq:circle-main-spectral-coefficients}.

\medskip
\noindent
\emph{Classification of the measures and polynomials.}
We apply the classification results for the four choices of
$\theta$.

\smallskip
\noindent
For $\theta=1$, Theorem~\ref{thm:circle-H-exclusion}
excludes the existence of such a measure.

For $\theta=z$, Theorem~\ref{thm:classification} gives
precisely the Lebesgue and Poisson families, with monic orthogonal
polynomials
\[
    \Phi_n(z)=z^n
    \qquad\text{or}\qquad
    \Phi_n(z)=z^{n-1}(z-a),\quad 0<|a|<1,
\]
respectively.

For $\theta=z(1-z)$, Corollary~\ref{cor:circle-G-complete}
gives the complete classification. More precisely, the
finite Verblunsky case is described by
Theorem~\ref{thm:circle-finite-G}, which gives the Poisson
family, including the normalized Lebesgue measure as $a=0$. The infinite
Verblunsky case is described by
Theorem~\ref{thm:circle-infinite-G}, which gives $\nu_q$
for a fixed $q\ne0$ satisfying $\operatorname{Re}q>-1/2$.
Theorem~\ref{thm:circle-Jacobi-branch} identifies its
monic orthogonal polynomials as $P_n^{(q)}$.

For $\theta=z^2$, Theorem~\ref{thm:circle-E-classification}
gives either the normalized Lebesgue measure or $\mu_{m,a}$, with $m\ge1$ and
$0<|a|<1$. In the latter case, the monic orthogonal
polynomials are
\[
    \Pi_n^{(m,a)}(z)=
    \begin{cases}
        z^n,&0\le n<m,\\
        z^{n-m}(z^m-a),&n\ge m.
    \end{cases}
\]
These results exhaust the possibilities, including the
solvable resonant cases.

\medskip
\noindent
\emph{Verification and characterization of the differential
coefficients.}
We verify that every listed family satisfies
\eqref{eq:circle-main-DE} and that the coefficient choices in
the table are exhaustive.

For $\Phi_n(z)=z^n$ and $n\ge2$, direct substitution gives
\[
    c_n=
    \begin{cases}
        1-n,&\theta=z\text{ or }\theta=z(1-z),\\
        0,&\theta=z^2,
    \end{cases}
\]
with $b_n$ arbitrary. This also covers the degrees $n<m$
of the family $\mu_{m,a}$.

For $\Phi_n(z)=z^{n-1}(z-a)$, write
$\theta(z)=z(1-\sigma z)$, where $\sigma=0$ or $1$.
Substitution into \eqref{eq:circle-main-DE}, using
\eqref{eq:circle-main-spectral-coefficients}, gives
for $n\ge2$
\[
    n-2+c_n=0,\qquad
    n(n-1+c_n)-a\bigl(b_n+2\sigma(n-1)\bigr)=0.
\]
Consequently,
\[
    c_n=2-n,\qquad b_n=\frac{n}{a}-2\sigma(n-1),
\]
as displayed in the table.

For the circular hypergeometric family, the differential
equation in Theorem~\ref{thm:circle-Jacobi-branch} yields
\[
    b_n=q+2-n,\qquad c_n=1-n-\overline q.
\]

For $\theta=z^2$ and
$\Phi_n(z)=z^n-a z^{n-m}$ with $n\ge m$ and $n\ge2$,
take $c_n=0$. By the spectral-coefficient formula,
\[
    \bigl(z^2D^2-b_nzD+\lambda_n\bigr)z^j
    =(j-n)(j+n-1-b_n)z^j,
    \qquad D=\frac{d}{dz}.
\]
Thus both powers $z^n$ and $z^{n-m}$ satisfy the equation
when $b_n=2n-m-1$. This gives the asserted coefficients.

To establish exhaustiveness, suppose that the same
$\Phi_n$ satisfies \eqref{eq:circle-main-DE} for two pairs
$(b_n,c_n)$ and $(\widetilde b_n,\widetilde c_n)$.
By \eqref{eq:circle-main-spectral-coefficients}, subtraction
of the two equations yields
\[
    (\Delta c-\Delta b\,z)\Phi_n'
         +n\Delta b\,\Phi_n=0,
\]
where $\Delta b=\widetilde b_n-b_n$ and
$\Delta c=\widetilde c_n-c_n$.
If $\Delta b=0$, then $\Delta c=0$.
Otherwise,
\[
    \left(z-\frac{\Delta c}{\Delta b}\right)\Phi_n'
       =n\Phi_n,
\]
whose monic polynomial solution is
\[
    \Phi_n(z)=\left(z-\frac{\Delta c}{\Delta b}\right)^n.
\]
Hence the differential coefficients are unique whenever
$\Phi_n$ has at least two distinct zeros.

For $n\ge2$, every nonmonomial polynomial in the table has
at least two distinct zeros. For the circular hypergeometric
family, one has
\[
    P_n^{(q)}(0)=\frac{(\overline q)_n}{(q+1)_n}\ne0.
\]
Moreover, all of its zeros lie in $\mathbb D$ and are
ordinary points of the Gauss differential equation, since
they are distinct from $0$ and $1$. Uniqueness for solutions
of a second-order linear differential equation therefore
implies that its zeros are simple.
The displayed coefficients for all the nonmonomial
polynomials are consequently unique.

Finally, at degree one, substitution of
$\Phi_1(z)=z+\Phi_1(0)$ gives
\[
    \lambda_1=b_1,\qquad c_1=-b_1\Phi_1(0),
\]
with $b_1$ arbitrary, independently of $\theta$.
This proves \eqref{eq:circle-main-degree-one} and completes
the verification of every admissible row in all positive
degrees.

The explicit densities of $\mu_{\mathrm{Leb}}$,
$\mu_{m,a}$, and $\nu_q$ establish the asserted absolute
continuity.
\end{proof}

\section{Concluding remarks: Pearson equations and signed measures}
\label{sec:circle-signed-measures}

 \paragraph{The Pearson equation approach.}
In the classical theory of orthogonal polynomials on the real line,
Pearson equations constitute one of the principal tools for
classification problems involving differential operators. In particular,
the existence of a fixed second-order differential operator
$\theta D^2+\psi D$ having the orthogonal polynomials as eigenfunctions
leads to the distributional equation
\[
    D(\theta u)=\psi u,
\]
where $u$ denotes the corresponding orthogonality functional.

The situation on the unit circle is considerably different.
Apart from the fact that orthogonality is formulated in terms of
Laurent polynomials, the coefficients of the differential equations
considered in the present work may depend on the degree $n$.
To illustrate this difficulty, consider
\[
    L_n=\theta(z)D_z^2+\psi_n(z)D_z,
    \qquad L_n\Phi_n=-\lambda_n\Phi_n,
\]
and let $u$ be the Hermitian orthogonality functional. Setting
$u_n=z^{1-n}u$, the orthogonality conditions become
\[
    \langle u_n,\Phi_n R\rangle=0,
    \qquad \deg R<n.
\]
Since $L_n$ preserves the degree filtration, the differential equation
and orthogonality imply
\[
    \left\langle
    D_z(\theta u_n)-\psi_n u_n,\,
    \Phi_n'R-\Phi_n R'
    \right\rangle=0,
    \qquad \deg R<n.
\]
However, the passage from these restricted identities to a full
Pearson equation
\[
    D_z(\theta u_n)=\psi_n u_n
\]
requires an additional argument, particularly because both the
functional $u_n$ and the admissible test space depend on $n$.

A noteworthy illustration of this difficulty is provided by
Theorem~\ref{thm:circle-infinite-G}. In the infinite branch of the
Gauss case, where $\theta(z)=z(1-z)$ and
$\psi_n(z)=c_n-b_nz$, the formal Pearson equation, expressed
in terms of the orthogonality measure, takes the form
\[
    D_\varphi[(1-e^{i\varphi})\mu]
    =i\bigl[(n-1+c_n)-(n-1+b_n)e^{i\varphi}\bigr]\mu.
\]
Remarkably, Theorem~\ref{thm:circle-infinite-G} establishes that
\[
    n-1+c_n=-\overline q,
    \qquad
    n-1+b_n=q+1,
    \qquad n\ge2,
\]
so that the dependence on $n$ disappears entirely and the Pearson
equation reduces to \eqref{eq:circle-G-Pearson}.

It is important to emphasize, however, that this rigidity is obtained
\emph{a posteriori}, only after the orthogonality measure has been
identified through the alternative and considerably involved
classification argument developed in the preceding sections.
Thus, the existence of a common Pearson equation emerges as a
consequence of the classification rather than as its starting point.

Despite its effectiveness in the classical real-line setting,
the Pearson equation approach does not therefore appear to provide
an immediate route to the classification problem considered here.
The longstanding difficulty of obtaining a complete classification
on the unit circle further illustrates the limitations of a direct
transfer of the classical techniques. For these reasons, we have
adopted the alternative approach developed in the present work.

This observation does not preclude the possibility of establishing
the classification by means of suitably adapted Pearson equations.
Indeed, determining whether a degree-independent Pearson-type
characterization can be deduced directly from the degree-dependent
differential equations for OPUC constitutes a problem of independent
interest and a natural direction for further research.

\paragraph{Extension to signed measures.}
The classification obtained in the preceding sections concerns positive
measures on the unit circle. A natural extension is to consider finite
real signed Borel measures $\mu$ on $\mathbb T$. In this setting,
the appropriate nondegeneracy condition is quasi-definiteness, namely,
\[
    \det T_N\ne 0,\qquad N\ge 0,
    \qquad T_N=(c_{j-k})_{j,k=0}^{N},
    \qquad c_k=\int_{\mathbb T}z^k\,d\mu(z),
\]
which ensures the existence of a unique sequence of monic orthogonal
polynomials $\{\Phi_n\}_{n\ge 0}$.

Although the Szeg\H{o} recurrence relation
\[
    \Phi_{n+1}(z)=z\Phi_n(z)+\delta_{n+1}\Phi_n^*(z),
    \qquad \delta_n=\Phi_n(0),
\]
remains valid, the restriction on the Verblunsky coefficients becomes
$|\alpha_n|\ne 1$, rather than $|\alpha_n|<1$. Indeed, setting
\[
    h_n=\int_{\mathbb T}|\Phi_n(z)|^2\,d\mu(z)\ne 0,
    \qquad \alpha_n=-\overline{\delta_{n+1}},
\]
one still has
\[
    h_{n+1}=(1-|\alpha_n|^2)h_n.
\]
Thus the Verblunsky coefficients may lie outside the unit disk, and the
zeros of the corresponding orthogonal polynomials need no longer be
confined to $\mathbb D$.

Consequently, the techniques developed in this work for
\[
    \theta(z)\in\{1,z,z^2,z(1-z)\}
\]
do not extend directly to the signed-measure setting. In particular,
arguments relying on the strict modulus bound for the Verblunsky
coefficients or on the localization of the zeros in $\mathbb D$ require
substantial modification. The classification for signed measures calls
for a separate treatment, taking into account both the indefinite
orthogonality conditions and the existence of finite signed representing
measures. This problem will be addressed in a subsequent publication.

 \section*{Declaration of the use of generative AI}

During the preparation of this manuscript, the author used
ChatGPT (OpenAI) as an auxiliary tool for the exploration
and refinement of mathematical arguments, the checking of
calculations, bibliographic searches, and the improvement
of the exposition. All AI-generated suggestions were
critically evaluated and independently verified by the
author, who takes full responsibility for the originality,
correctness, and integrity of the work.

  \section*{Acknowledgments }

The author wishes to  thank  Professor Francisco Marcell\'an for his valuable comments. 

\bibliography{opuc}


\end{document}